\documentclass[leqno,12pt]{amsart}
\usepackage{amssymb}
\usepackage{amsfonts,amsthm}
\usepackage{hyperref}

\newtheorem{theorem}{Theorem}[section]
\newtheorem{corollary}[theorem]{Corollary}
\newtheorem{lemma}[theorem]{Lemma}
\newtheorem{proposition}[theorem]{Proposition}
\theoremstyle{definition}
\newtheorem{defn}[theorem]{Definition}
\theoremstyle{remark}
\newtheorem*{remark}{Remark}
\newtheorem{exmp}{Example}

\newcommand{\nc}{\newcommand}
\newcommand{\be}{\begin{equation}}
\newcommand{\ee}{\end{equation}}
\newcommand{\no}{\nonumber}
\newcommand{\bc}{\begin{center}}
\newcommand{\ec}{\end{center}}
\nc{\bth}{\begin{theorem}}
\nc{\bpr}{\begin{proposition}}
\nc{\epr}{\end{proposition}}
\nc{\ble}{\begin{lemma}}
\nc{\ele}{\end{lemma}}
\nc{\bco}{\begin{corollary}}
\nc{\eco}{\end{corollary}}
\nc{\bre}{\begin{remark}}
\nc{\ere}{\end{remark}}
\nc{\f}{\frac}
\nc{\na}{\nabla}
\nc{\al}{\alpha}
\nc{\bet}{\beta}
\nc{\ga}{\gamma}
\nc{\de}{\delta}
\nc{\ep}{\varepsilon}
\nc{\ei}{\ep _i}
\nc{\ea}{\ep _{\al}}
\nc{\ej}{\ep _j}
\begin{document}
\title[Semi-Riemannian submersions with totally geodesic fibres]
{      Semi-Riemannian submersions with\\
       totally geodesic fibres}
\author{Gabriel B\u adi\c toiu}
\address{
    Institute of Mathematics
    of the Romanian Academy\newline
 \indent    P.O. Box 1-764\newline
 \indent    Bucharest 014700\newline
 \indent    Romania}
\maketitle
\begin{abstract}
        We classify semi-Riemannian submersions with connected totally
        geodesic fibres
        from a real pseudo-hyperbolic space onto a semi-Riemannian
        manifold under the assumption that the dimension of the fibres
        is less than or equal to three and the metrics induced on fibres are negative definite.
         Also, we obtain the
        classification of semi-Riemannian
        submersions with connected complex
        totally geodesic fibres from a complex pseudo-hyperbolic space onto a
        semi-Riemannian manifold under the assumption that the dimension of the fibres
        is less than or equal to two and the metric induced on fibres are negative definite.
        We prove that there are no semi-Riemannian submersions
        with connected quaternionic fibres
        from a quaternionic pseudo-hyperbolic space onto a Riemannian manifold.
\end{abstract}
\section{Introduction and main results}

Riemannian submersions, introduced
by O'Neill \cite{one} and Gray \cite{gra},
 have been used by many authors to
construct specific Riemannian metrics.
A systematic exposition can be found
in Besse's book \cite{bes}.
In this paper, we obtain classification results for semi-Riemannian
submersions with totally geodesic fibres.

We first recall briefly some related work on the classification
problem of semi-Riemannian submersions.
Escobales \cite{esc, esco} and
Ranjan  \cite{ran} classified Riemannian submersions
with connected totally geodesic fibres from an $n$-sphere $S^n$, and
with  connected complex  totally geodesic fibres from a complex projective $n$-space
$\mathbb CP^n$, respectively. Ucci \cite{ucc} showed  that there are no Riemannian submersions
with fibres $\mathbb CP^3$ from the complex projective space $\mathbb CP^7$ onto $S^8(4)$,
and with fibres $\mathbb HP^1$  from the
quaternionic projective space $\mathbb HP^3$ onto $S^8(4)$.
    In \cite{ranj}, Ranjan obtained a classification theorem for  Riemannian
submersions with connected totally geodesic fibres from a compact simple Lie group.
Gromoll and Grove obtained  in \cite{gro} that, up to equivalence, the only
Riemannian submersions of spheres (with connected fibres)
are the Hopf fibrations, except possibly
for fibrations of the $15$-sphere by homotopy $7$-spheres.
This classification was invoked in the
proof of the Diameter Rigidity Theorem (see \cite{grom})
and of the Radius Rigidity Theorem
(see \cite{wil}).
Using an approach different from Gromoll and Grove \cite{gro}, Wilking
\cite{wilk} proved that a Riemannian submersion $\pi:S^m\to B^b$ is
metrically equivalent to the Hopf fibration for $(m,b)=(15,8)$ and
obtained an improved version of the Diameter Rigidity Theorem
as a consequence of his classification theorem.

In comparison, there are few classification results for semi-Riemannian
submersions, and the consequences seem to be at least as important
as those for Riemannian submersions.
In \cite{mag}, Magid proved that the only semi-Riemannian submersions with totally
geodesic fibres from an anti-de Sitter space onto
 a Riemannian manifold are the canonical
semi-Riemannian submersions $H^{2m+1}_1\to\mathbb CH^m$.
In \cite{bad}, the present author and Stere Ianu\c s classified  semi-Riemannian submersions
with connected totally geodesic fibres from a pseudo-hyperbolic space onto a Riemannian manifold,
and with  connected complex totally geodesic fibres from a complex pseudo-hyperbolic space onto
a Riemannian manifold.

The aim of this work is to prove new classification
results in the theory of semi-Riemannian submersions analogous to
those in Riemannian geometry.
It is my pleasure to thank Professor Stere Ianu\c s
for useful discussions on this subject.

Now, we list the main results proved in this paper.

\bth\label{t:1.1}
    Let $\pi:H^{n+r}_{s+r'}\to B^n_s$
        be a semi-Riemannian submersion with connected totally
    geodesic fibres from a pseudo-hyperbolic space
        onto a semi-Riemannian manifold.
    If the dimension of the fibres is less than or equal
        to $3$ and if the metrics induced on fibres are negative definite, then $\pi$ is equivalent to one of
    the following canonical semi-Riemannian submersions$:$
    \begin{itemize}
    \item[(a)] $H^{2m+1}_{2t+1}\to\mathbb CH^m_t,\ \  0\leq t\leq m.$
    \item[(b)] $H^{4m+3}_{4t+3}\to\mathbb HH^m_t,\ \  0\leq t\leq m.$
\end{itemize}
\end{theorem}
\bth\label{t:1.2}
    Let $\pi:H^{n+r}_{s+r'}\to B^n_s$ be a semi-Riemannian submersion with connected totally
    geodesic fibres from a pseudo-hyperbolic space onto a semi-Riemannian manifold.
    Assume that one of the following conditions is satisfied$:$
    \begin{itemize}
    \item[(A)] $B$ is an isotropic semi-Riemannian manifold, which means that
      for any
      $x\in B_s^n$ and any real number $t$, the group of isometries
      $\mathbf I(B^n_s,g')$
       preserving $x$ acts transitively on the set of
       all nonzero tangent vectors $X$ at $x$ for which $g'(X,X)=t$, or
    \item[(B)] $\mathrm{index\ }(B)\in\{0,\dim B\}$.
    \end{itemize}
    Then $\pi$ is equivalent to one of the following
    canonical semi-Riemannian submersions$:$
    \begin{itemize}
    \item[(a)] $H^{2m+1}_{2t+1}\to\mathbb CH^m_t,\ \ 0\leq t\leq m$.
    \item[(b)] $H^{4m+3}_{4t+3}\to\mathbb HH^m_t,\ \ 0\leq t\leq m$.
    \item[(c)] $H^{15}_{7+8t}\to H^8_{8t}(-4),\ \ t\in\{0,1\}$.
    \end{itemize}
\end{theorem}
\bth\label{t:1.3}
    Let $\pi:\mathbb CH^{n}_{s}\to B$ be a semi-Riemannian submersion
    from a complex pseudo-hyperbolic space onto a semi-Riemannian manifold.
    Assume that the fibres are connected  complex
    totally geodesic submanifolds,
    and one of the following conditions  is satisfied$:$
    \begin{itemize}
        \item[(A)] The real dimension of the fibres is $r\leq 2$ and the fibres are negative definite, or
        \item[(B)] $B$ is an isotropic semi-Riemannian manifold, or
        \item[(C)] $\mathrm{index\, }(B)\in\{0,\dim B\}$.
    \end{itemize}
    Then $\pi$ is equivalent to the canonical semi-Riemannian submersion
    \begin{itemize}
        \item[] $\mathbb CH^{2m+1}_{2t+1}\to\mathbb HH^m_t,\ \ 0\leq t\leq m$.
    \end{itemize}
\end{theorem}
\bth\label{t:1.4}
There exist no semi-Riemannian submersions  $\pi:\mathbb HH^n_s\to B$
with connected quaternionic fibres
from a quaternionic pseudo-hyperbolic space onto an isotropic semi-Riemannian manifold
or onto a semi-Riemannian manifold of $\mathrm{index}(B)\in\{0,\dim(B)\}$.
\end{theorem}

\section{Preliminaries and examples}
In this section we recall several notions
and results which will be needed throughout the paper.
We also exhibit the construction of canonical semi-Riemannian submersions.
\begin{defn}
    Let $(M,g)$ be an $(n+r)$-dimensional connected
    semi-Riemannian manifold of index $s+r'$, and
    $(B,g')$ an $n$-dimensional connected
    semi-Riemannian manifold of index $s$, where $0 \leq s\leq n$, $0 \leq r'\leq r$.
    A {\it semi-Riemannian submersion} (see \cite{onei}) is
    a smooth map
         $\pi :M\to B$
    which is surjective and satisfies the following axioms:
\begin{itemize}
\item[(a)] $\pi_*|_p$ is surjective for all $p\in M$;
\item[(b)] the fibres $\pi^{-1}(b)\ , \ b\in B$, are
    semi-Riemannian submanifolds of $M$;
\item[(c)] $\pi_*$ preserves scalar products of
    vectors normal to fibres.
\end{itemize}
\end{defn}

\vspace*{7pt}
  We shall always assume that the fibres are connected, the dimension of the fibres
  $\dim M-\dim B>0$  and $\dim B>0$.
   The vectors tangent to fibres are called
    vertical and those normal to fibres
    are called horizontal.
    We denote by $\mathcal V$  the vertical distribution and by $\mathcal H$
    the horizontal distribution.

    The geometry of  semi-Riemannian submersions is characterized by
    O'Neill's  tensors $T$, $A$ (see \cite{one}, \cite{onei})
    defined for vector fields  $E$, $F$ on $M$ by
\begin{eqnarray*}
        A_EF &=&h \na _{hE}{vF} + v \na _{hE}{hF},\\
        T_EF &=&h \na _{vE}{vF} + v \na _{vE}{hF},
\end{eqnarray*}
    where $\na $ is the Levi-Civita connection of $g$, and $v$ and $h$ denote the orthogonal
    projections on $\mathcal V$  and $\mathcal H$, respectively.
    For basic properties of O'Neill's tensors see \cite{one}, \cite{onei}, \cite{bes} or
    \cite{ian}.

\begin{defn}
    (i) A vector field $X$ on $M$ is said to be {\it basic}
    if $X$ is horizontal and $\pi$-related to
    a vector field $X'$ on $B$.

    (ii) A vector field $X$ along the fibre
    $\pi^{-1}(x)$, $x\in B$, is said to be {\it basic along} $\pi^{-1}(x)$
    if  $X$ is horizontal and $\pi_{*p}X(p)=\pi_{*q}X(q)$
    for every $p$, $q\in\pi^{-1}(x)$.
\end{defn}
We notice that each vector field $X'$
on $B$ has a unique horizontal lift $X$ to $M$
which is basic. For a vertical vector field $V$ and  a  basic vector field $X$
we have $h\na_VX=A_XV$ (see \cite{one}).
We denote by  $R$, $R'$ and  $\hat R$
the Riemann curvature tensors of $M$, $B$ and
of the fibre $\pi^{-1}(x)$, $x\in M$, respectively.
We choose the convention for the curvature tensor
$R(E,F)=\na_E\na_F-\na_F\na_E-\na_{[E,F]}$. The Riemann curvature tensor is defined by
$$R(E,F,G,H)=g(R(G,H)F,E).$$
For O'Neill's equations of a semi-Riemannian submersion we refer to \cite{one} or
\cite{bes}.

\begin{defn}
     Two semi-Riemannian submersions $\pi,\pi':(M,g)\to (B,g')$
     are said to be {\it equivalent} if there exists an
     isometry $f$ of $M$ which induces an isometry $\tilde f$ of $B$ so that
     $\pi'\circ f=\tilde f\circ\pi$. The pair $(f,\tilde f)$ is called a
     {\it bundle isometry}.
\end{defn}

We shall need the following theorem, which is the semi-Riemannian version of
Theorem 2.2 in \cite{esc}.
\bth\label{t:4}
     Let $\pi_1,\pi_2:M\to B$ be semi-Riemannian submersions from a
     complete connected semi-Riemannian manifold $M$ onto a semi-Riemannian manifold $B$.
     Assume that the fibres of these submersions are connected and totally
     geodesic and the metric induced on fibres is negative definite.
     Let $f$ be an isometry of $M$ satisfying the following
     properties at a given point $p\in M:$
\begin{itemize}
\item[(1)] $f_{*p}:T_pM\to T_{f(p)}M$ maps $\mathcal{H}_{1p}$ onto
      $\mathcal{H}_{2f(p)}$, where  $\mathcal{H}_{i}$ denote the horizontal
     distributions of $\pi_i$ for  $i\in\{1,2\}$.
\item[(2)] $f_*A_{1E}F=A_{2{f_*E}}f_*F$ for every $E$, $F\in T_pM$,
     where $A_i$ are the integrability tensors associated with $\pi_i$.
\end{itemize}
     Then $f$ induces an isometry $\Tilde{f}$ of $B$ so that the pair
     $(f,\Tilde{f})$ is a bundle isometry between $\pi_1$ and $\pi_2$.
     In particular, $\pi_1$ and $\pi_2$ are equivalent.
\end{theorem}

Escobales's proof of Theorem 2.2 in \cite{esc},  also works in this
semi-Riemannian case. He proves that for any $b\in B$
 which can be joined with $\pi_1(p)$
by a geodesic we have:
\begin{itemize}
  \item[(i)] for every $x\in\pi_1^{-1}(b)$, $f_{*x}:T_xM\to T_{f(x)}M$ maps
     $\mathcal H_{1x}$ onto $\mathcal H_{2f(x)}$, and
  \item[(ii)] $f$ maps the fibre $\pi_1^{-1}(b)$ into the fibre $\pi_2^{-1}(\pi_2(f(x)))$ with
    $x\in\pi_1^{-1}(b)$.
\end{itemize}
We  notice that for any $x\in\pi_1^{-1}(b)$ with $b\in B$,
which can be joined with $\pi_1(p)$ by a geodesic, the conditions (1) and (2) are also satisfied
for the point $x$.
Since $M$ is connected, $B$ is also connected. Therefore, any point $b\in B$
can be joined with $\pi_1(p)$ by a broken geodesic. Repeating the argument above, for any corner point
of this broken geodesic, we see that for any $b\in B$, $f$ maps the fibre $\pi_1^{-1}(b)$
into a fibre.

\vspace*{7pt}
\begin{defn}
   Let $\langle\cdot,\cdot\rangle$ be the symmetric bilinear form on $\mathbb R^{m+1}$
   given by
        $$\langle x,y\rangle=-\sum\limits_{i=0}^sx_iy_i+\sum\limits_{i=s+1}^mx_iy_i$$
  for $x=(x_0,\cdot\cdot\cdot,x_m),y=(y_0,\cdot\cdot\cdot,y_m)\in\mathbb R^{m+1}$.
  For any $c<0$ and any positive integer $s$, let
  $H^m_s(c)=\{x\in\mathbb R^{m+1}\ |\ \langle x,x\rangle =1/c\}$ be the semi-Riemannian submanifold of
    $$\mathbb R_{s+1}^{m+1}=(\mathbb R^{m+1},
  ds^2=-dx^0\otimes dx^0-\cdot\cdot\cdot -dx^s\otimes dx^s+
        dx^{s+1}\otimes dx^{s+1}+\cdot\cdot\cdot +dx^m\otimes dx^m).$$
  $H^m_s(c)$ is called the $m$-dimensional
  ({\it real}) {\it pseudo-hyperbolic space} of index $s$.
\end{defn}

\vspace*{7pt}
  We notice that $H^m_s(c)$ has constant sectional curvature $c$,
  whose curvature tensor is given by
  $R(X,Y,X,Y)=c(g(X,X)g(Y,Y)-g(X,Y)^2)$.
  We shall denote simply  $H^m_s=H^m_s(-1)$.
  It should be remarked that $H^m_s$ can be written as a homogeneous space,
  namely
  $H^m_s=SO(s+1,m-s)/SO(s,m-s)$, $H^{2m+1}_{2s+1}=SU(s+1,m-s)/SU(s,m-s)$, and
  $H^{4m+3}_{4s+3}=Sp(s+1,m-s)/Sp(s,m-s)$ (see \cite{wol}).

\vspace*{7pt}
\begin{defn}\label{d:3}
   Let $(\cdot,\cdot)$ be the Hermitian form on $\mathbb C^{m+1}$ given by
     $$(z,w)=-\sum\limits_{i=0}^sz_i\bar{w_i}+\sum\limits_{i=s+1}^mz_i\bar{w_i}$$
  for $z=(z_0,\cdot\cdot\cdot,z_m), w=(w_0,\cdot\cdot\cdot,w_m)\in\mathbb C^{m+1}$.
  For $c<0$, let $M(c)$ be the real hypersurface of $\mathbb C^{m+1}$
  given by $M(c)=\{z\in\mathbb C^{m+1}\ |\ (z,z)= 4/c\}$, which is
  endowed with the induced metric of
  $$(\mathbb C^{m+1},
  ds^2=-dz^0\otimes d\bar{z}^0-\cdot\cdot\cdot -dz^s\otimes d\bar{z}^s+
        dz^{s+1}\otimes d\bar{z}^{s+1}+\cdot\cdot\cdot +dz^m\otimes d\bar{z}^m).$$
  The natural action of $S^1=\{e^{i\theta}\ |\ \theta\in\mathbb R\}$ on $\mathbb C^{m+1}$
  induces an action on $M(c)$. Let $\mathbb CH^m_s(c)=M(c)/S^1$
  endowed with the unique indefinite K\"ahler
  metric of index $2s$ such that the projection $M(c)\to M(c)/S^1$
  becomes a semi-Riemannian submersion
  (see \cite{bar}).
  $\mathbb CH^m_s(c)$ is called the {\it complex pseudo-hyperbolic space}.
\end{defn}

\vspace*{7pt}
  Notice that $\mathbb CH^m_s(c)$ has constant holomorphic sectional curvature $c$,
  whose curvature tensor is given by
  $R(X,Y,X,Y)= (c/4)(g(X,X)g(Y,Y)-g(X,Y)^2+3g(I_0X,Y)^2)$, where
  $I_0$ is the natural complex structure on $\mathbb CH^m_s(c)$.
  We shall denote simply
   $\mathbb CH^m_s=\mathbb CH^m_s(-4)$. It is well-known that
   $\mathbb CH^m_s$ is a homogeneous space, namely
   $\mathbb CH^m_s=SU(s+1,m-s)/S(U(1)U(s,m-s))$ and
   $\mathbb CH^{2m+1}_{2s+1}=Sp(s+1,m-s)/U(1)Sp(s,m-s)$  (see \cite{wol}).

We shall denote by $\mathbb HH^n_s$ the quaternionic pseudo-hyperbolic
space of real dimension $4n$, and of
quaternionic index $s$ with quaternionic sectional curvature $-4$, and by
$S^n$ and $S^n(4)$ the spheres with sectional curvature $1$ and $4$, respectively.

By a standard construction (see Theorem 9.80 in \cite{bes}),
one can obtain many examples of
semi-Riemannian submersions with totally geodesic fibres of type
$\pi:G/K\to G/H$,  where $G$ is a Lie group and $K$, $H$
 are
closed Lie subgroups of $G$ with $K\subset H$.
In this way the following {\it canonical semi-Riemannian submersions}, also
called  {\it generalized Hopf fibrations}, are obtained:\\

\begin{exmp}
    Let $G=SU(t+1,m-t)$, $H=S(U(1)U(t,m-t))$, $K=SU(t,m-t)$.
    For every $0\leq t\leq m$, we have the semi-Riemannian submersion
    $$H^{2m+1}_{2t+1}=SU(t+1,m-t)/SU(t,m-t)\to
       \mathbb CH^m_t=SU(t+1,m-t)/S(U(1)U(t,m-t)).$$
\end{exmp}

\begin{exmp}
    Let $G=Sp(t+1,m-t)$, $H=Sp(1)Sp(t,m-t)$, $K=Sp(t,m-t)$.
For every $0\leq t\leq m$, we get the semi-Riemannian submersion
   $$H^{4m+3}_{4t+3}=Sp(t+1,m-t)/Sp(t,m-t)\to
       \mathbb HH^m_t=Sp(t+1,m-t)/Sp(1)Sp(t,m-t).$$
\end{exmp}

\begin{exmp}
   a)\ Let $G=Spin(1,8)$, $H=Spin(8)$, $K=Spin(7)$. Then
   we have the semi-Riemannian submersion (see \cite{bad})
   $$H^{15}_7=Spin(1,8)/Spin(7)\to H^8(-4)=Spin(1,8)/Spin(8).$$
   b)\ Let $G=Spin(9)$, $H=Spin(8)$, $K=Spin(7)$. Then
   we have the semi-Riemannian submersion (see \cite{bes})
   $$S^{15}=Spin(9)/Spin(7)\to S^8(4)=Spin(9)/Spin(8).$$
\end{exmp}

\begin{exmp}
       Let $G=Sp(t+1,m-t)$, $H=Sp(1)Sp(t,m-t)$, $K=U(1)Sp(t,m-t)$.
    For every $0\leq t\leq m$, we obtain the semi-Riemannian
    submersion
    $$\mathbb CH^{2m+1}_{2t+1}=Sp(t+1,m-t)/U(1)Sp(t,m-t)\to
    \mathbb HH^m_t=Sp(t+1,m-t)/Sp(1)Sp(t,m-t).$$
\end{exmp}

\vspace*{7pt}
In order to prove Theorem \ref{t:1.2}, we need the following nonexistence proposition,
which is the semi-Riemannian version of Proposition 5.1 in \cite{ran}.
\bpr\label{cayley}
    There exist no semi-Riemannian submersions $\pi:H^{23}_{7+8t}\to\mathbb CaH^2_t$,
    $t\in\{0,1,2\}$, with totally geodesic fibres from the $23$-dimensional pseudo-hyperbolic
    space of index $7+8t$ onto the Cayley pseudo-hyperbolic plane of Cayley index $t$ .
\epr
We notice that the case $t=2$ is Proposition 5.1 in \cite{ran}.
For the case $t=0$, see \cite{bad}.
Here we only recall  some details of Ranjan's proof and suggest its modification
to the semi-Riemannian case.
Ranjan's argument in \cite{ran}, which leads to a
contradiction to the assumption of the existence
of such a submersion, is based on finding for every
$X\in\mathcal H_p$, $g(X,X)\not=0$,
an irreducible $Cl(\mathcal V_p)$-submodule $S$ of
$\mathcal H_p$ passing through $X$.
Here $Cl(\mathcal V_p)$ denotes the Clifford algebra
of $(\mathcal V_p,\Tilde g_p)$, where $\Tilde g(U,V)=-g(U,V)$
for every $U$, $V\in\mathcal V_p$.
$\mathcal H_p$ becomes a $Cl(\mathcal V_p)$-module
by considering the extension
of the map
$\mathcal U:\mathcal V_p\to\mathrm{End}(\mathcal H_p)$
defined by
 $\mathcal U(V)(X)=A_XV$
to the Clifford algebra $Cl(\mathcal V_p)$.
Since $\Tilde g_p$ is positive definite, we have
$Cl(\mathcal V_p)\simeq\mathbb R(8)\oplus\mathbb R(8)$.
Hence, $\mathcal H_p$ splits into
two $8$-dimensional irreducible $Cl(\mathcal V_p)$-modules.
Since the induced metrics on
fibres are negative definite,
we obtain in a manner similar
to Ranjan's proof that
\begin{itemize}
     \item[(i)\ ] for $g(X,X)>0$, $\pi^{-1}(\mathbb CaH^1)$
       is totally geodesic in $H^{23}_{7+8t}$ and is isometric
       to $H^{15}_7$, where $\mathbb CaH^1$ denotes the
       Cayley hyperbolic line through $\pi_*X$, and
     \item[(ii)\ ]  for $g(X,X)<0$, $\pi^{-1}(\mathbb CaH^1_1)$
       is totally geodesic in $H^{23}_{7+8t}$ and is isometric
       to $H^{15}_{15}$, where $\mathbb CaH^1_1$ denotes the
       negative definite Cayley hyperbolic line through $\pi_*X$.
\end{itemize}
We choose $S$ to be the horizontal space of the restricted submersion
$\Tilde\pi:H^{15}_7\to\mathbb CaH^1$ if $g(X,X)>0$ or
$\Tilde\pi:H^{15}_{15}\to\mathbb CaH^1_1$ if $g(X,X)<0$.

\section{Proof of the main results}

The next lemma gives useful properties of
O'Neill's integrability tensor.
\ble\label{lem:2}
    Let $\pi:M\to B$ be a semi-Riemannian submersion
    with connected totally geodesic
    fibres from a semi-Riemannian manifold $M$ with
    constant curvature
    $c\neq 0$. Then the following assertions are true$:$
\begin{itemize}
    \item[(a)]
    If $X$ is a horizontal vector such that $g(X,X)\neq 0$, then
    the map $A_X:\mathcal V\to\mathcal H$
        given by $A_X(V)=A_XV$
        is injective and the map
    $A^*_X:\mathcal{H}\to\mathcal{V}$ given by
    $A^*_X(Y)=A_XY$ is surjective.
\item[(b)] If $X$, $Y$ are the horizontal liftings
    along the fibre
    $\pi^{-1}(\pi(p))$, $p\in M$, of
    two vectors $X', Y'\in T_{\pi(p)}B$ respectively, $g'(X',X')\neq 0$
    and $(A_XY)(p)=0$, then $A_XY=0$ along the fibre
    $\pi^{-1}(\pi(p))$.
\end{itemize}
\ele
\begin{proof}
(a)  By O'Neill's equations, we get
    $$g(A_XV, A_XW)=c g(X,X)g(V,W)$$
    for a horizontal vector field $X$ and for
    vertical vector fields $V$ and $W$.
    Thus
    $A^*_XA_XV=-c g(X,X)V$ for every  vertical vector field $V$.
    Therefore $A_X:\mathcal{V}\to\mathcal{H}$ is injective
    and $A^*_X:\mathcal{H}\to\mathcal{V}$ is surjective.

(b)  By O'Neill's equations, we have
    $$-3g(A_XY,A_XZ)=c[g(X,X)g(Y,Z)-g(X,Y)g(X,Z)]-
    R'(\pi_*X,\pi_*Y,\pi_*X,\pi_*Z)$$
    for horizontal vector fields $X$, $Y$ and $Z$.

    If $X$, $Y$, $Z$ are basic vector fields, then
    $g(A_XY, A_XZ)$ is constant along the fibre
    $\pi^{-1}(\pi(p))$. Therefore,
    $g(A_XA_XY,Z)=0$ along the fibre $\pi^{-1}(\pi(p))$
        for every  basic vector field $Z$. Hence $A_XA_XY=0$ along $\pi^{-1}(\pi(p))$.
    Since $A_{X}:\mathcal{V}\to\mathcal{H}$
        is injective, it follows that $A_XY=0$ along the fibre
    $\pi^{-1}(\pi(p))$.
\end{proof}

\ble\label{l:1}
     If $\pi:M\to B$ is a semi-Riemannian submersion with  connected totally geodesic
     fibres from a semi-Riemannian manifold $M$ with constant curvature
     $c\neq 0$ onto a semi-Riemannian manifold $B$, then the tangent bundle
     of any fibre is trivial.
\ele
\begin{proof}
      Let $x\in B$ and $p\in\pi^{-1}(x)$.
      Let
      $\{v_{1p},\dots  ,v_{rp}\}$ be an orthonormal basis in $\mathcal V_p$.
      Let $Y_1,Y_2,\dots,Y_r$ be the horizontal liftings along the fibre
      $\pi^{-1}(\pi(p))$ of
      $(1/(cg(X,X)))\pi_*A_Xv_{1p}$, $(1/(cg(X,X)))\pi_*A_Xv_{2p}$,\dots,
      $(1/(cg(X,X)))\pi_*A_Xv_{rp}$, respectively.
      Let $v_i=A_XY_i$ for each $i\in\{1,\dots,r\}$.
      Since
\begin{eqnarray*}
       g(v_j,v_l)&=&g(A_XY_j,A_XY_l)\\
       &=&
      (1/3)(R'(\pi_*X,\pi_*Y_j,\pi_*X,\pi_*Y_l)-cg(X,X)g(Y_j,Y_l)+cg(X,Y_j)g(X,Y_l))
\end{eqnarray*}
      is constant along the fibre
      $\pi^{-1}(\pi(p))$ and
       $$g(A_XY_j,A_XY_l)(p)=\f 1{c^2}g(A_XA_Xv_{jp},A_XA_Xv_{lp})=
      g(X,X)^2g(v_{jp},v_{lp})=\ej\de_{jl},$$
      we see that $\{v_1,v_2,\dots,v_r\}$ is a global orthonormal basis
      of the tangent bundle of the fibre $\pi^{-1}(x)$,
      which makes the tangent bundle trivial.
\end{proof}

We suppose that the curvature of the total space is negative.
The case of positive curvature can be reduced to the negative
one by changing simultaneously the signs of the metrics on the base
and on the total space.
We establish relations between the dimensions
and the indices of fibres and of base spaces, and see how the
geometry of base spaces looks like.

\bth\label{p:3}
    Let $\pi:M\to B$ be a semi-Riemannian submersion
    with  connected totally geodesic fibres from an $(n+r)$-dimensional
    semi-Riemannian manifold $M$ of index $s+r'$ with
    constant negative curvature $c$ onto an $n$-dimensional
    semi-Riemannian manifold $B$ of index $s$. Then the following hold$:$
\begin{itemize}
    \item[(1) ] $n=k(r+1)$ for some  positive integer $k$
    and $s=q_1(r'+1)+q_2(r-r')$ for
    some  nonnegative integers $q_1$, $q_2$ with
    $q_1+q_2=k$.
    \item[(2) ] If, moreover, $M$ is a simply connected complete
    semi-Riemannian manifold
    and the dimension of fibres is less than or equal to $3$
    and the metric induced on fibres is negative definite,
    then $B$ is an
    isotropic semi-Riemannian manifold and $r\in\{1,3\}$.
\end{itemize}
\end{theorem}
\begin{proof}
      Normalizing the metric on $M$, we can suppose $c=-1$.
      Let $p\in M$. Since the tangent bundle of the fibre
      $\pi^{-1}(\pi(p))$ is trivial, we can choose a global
      orthonormal frame $\{v_1,v_2,\dots ,v_r\}$
      for the tangent bundle of $\pi^{-1}(\pi(p))$. We have
      $g(v_i,v_j)=\ei\de_{ij}$,
      $\ei\in\{-1,1\}$, and card$\{i | \ei<0\}=r'$.

      (1) Let $X$ be the horizontal lifting along the fibre $\pi^{-1}(\pi(p))$
      of a vector $X'\in T_{\pi(p)}B$, so that $g(X',X')\in\{-1,1\}$.
      By O'Neill's equations, we have
         $$g(A_YV, A_YV)=-g(Y,Y)g(V,V)$$
      for a horizontal vector field $Y$ and for a vertical vector field $V$.
      Along the fibre $\pi^{-1}(\pi(p))$ we obtain for every $i,j\in\{1,\dots,r\}$
      $$g(A_Xv_i,A_Xv_j)=-g(X,X)g(v_i,v_j)=-g(X,X)\ei\de_{ij} , $$
      $$g(X,A_Xv_i)=-g(A_XX,v_i)=0.$$
      Thus $\{X,A_Xv_1,\dots  ,A_Xv_r\}$ is an orthonormal system. Hence $n\geq r+1$.

      Let $L_0=X$.
      For every  integer $\al$ such that $1\leq\al <n/(r+1)$,
      let $L_\al$ be a horizontal vector field along the fibre
      $\pi^{-1}(\pi(p))$ such that $L_\al$ is
      the horizontal lifting of some unit vector
      (i.e., $g(L_\al,L_\al)\in\{-1,1\}$), that $L_\al$ is orthogonal to
      $L_0,L_1,\dots  ,L_{\al-1}$ and that
      $L_\al(p)\in\ker A^*_{L_0(p)}\cap\ker A^*_{L_1(p)}\cap\dots
      \cap\ker A^*_{L_{\al-1}(p)}$. Then, by Lemma \ref{lem:2}, $L_\al(q)$ belongs to
      $\ker A^*_{L_0(q)}\cap\ker A^*_{L_1(q)}\cap\dots
      \cap\ker A^*_{L_{\al-1}(q)}$ for every $q\in\pi^{-1}(\pi(p))$.
      Therefore, for $j\in\{1,\dots,r\}$ and $\al,\beta\geq 0$,
      we get
      $$g(A_{L_\al}v_j,L_\beta)=-g(v_j,A_{L_\al}L_\beta)=0$$
      along the fibre $\pi^{-1}(\pi(p))$.

      By O'Neill's equations, we obtain
      \begin{equation}\label{e:ecu}\begin{split}
            R(X,U,Y,V)& =g((\na_UA)_XY,V)+g(A_XU,A_YV)\\
            & =g(\na_UA_XY,V)-g(A_{\na_UX}Y,V)-g(A_X\na_UY,V)+g(A_XU,A_YV)\\
            & =g(\na_UA_XY,V)+g(A_YA_XU,V)-g(A_XA_YU,V)-g(A_YA_XU,V)\\
            & =g(\na_UA_XY,V)+g(A_YU,A_XV)
      \end{split}\end{equation}
      for  basic vector fields $X$, $Y$
      and for vertical vector fields $U$, $V$.
      Thus, along the fibre $\pi^{-1}(\pi(p))$ we get for every
      $\al,\beta\geq 0$ and $j,l\in\{1,\dots,r\}$
      \begin{eqnarray*}
            g(A_{L_\al}v_j,A_{L_\beta}v_l)& =
                  &R(L_\al,v_l,L_\beta,v_j)-g(\na_{v_l}A_{L_\al}L_\beta,v_j)\\
            & =&-g(L_\al,L_\beta)g(v_l,v_j)-v_l(g(A_{L_\al}L_\beta,v_j))
               +g(A_{L_\al}L_\beta,\na_{v_l}v_j).
      \end{eqnarray*}
      Since $A_{L_\al}L_\beta=0$ along the fibre $\pi^{-1}(\pi(p))$, it follows that
      $$g(A_{L_\al}v_j,A_{L_\beta}v_l)=-g(L_\al,L_\beta)g(v_l,v_j)=
            -g(L_\al,L_\beta)\varepsilon_l\de_{lj}.$$
      We proved that for some positive integer  $k$,
        $$\mathcal{L}=\{L_0,A_{L_0}v_1,\dots,A_{L_0}v_r,\dots,
           L_{k-1}, A_{L_{k-1}}v_1,\dots  ,A_{L_{k-1}}v_r\}$$
      is an orthonormal basis of $\mathcal H$  along the fibre $\pi^{-1}(\pi(p))$.
      Thus $\dim B=(1+\dim\mathrm{ fibre})k$
      for some  positive integer  $k$.
      Counting the timelike vectors in $\mathcal{L}$, we get
      $\mathrm{index}(B)=q_1(r'+1)+q_2(r-r')$
      for some nonnegative integers  $q_1$, $q_2$ with $q_1+q_2=k$.

(2)\
      Let $x\in B$ and $X'$, $Y'\in T_xB$ such that $g'(X',X')=g'(Y',Y')\not=0$.
      We shall construct an isometry $\Tilde f:B\to B$ such that
      $\Tilde f(x)=x$ and $\Tilde f_*X'=Y'$.
      Note that we may assume that
      $g'(X',X')=g'(Y',Y')=\pm 1$.
      Let $X$, $Y$ be the horizontal liftings along the fibre
      $\pi^{-1}(x)$ of $X'$ and $Y'$, respectively.
      Take $p\in\pi^{-1}(x)$. Let
      $$\mathcal L=\{L_0,A_{L_0}v_1,\dots  ,A_{L_0}v_r,\dots ,
            L_{k-1},A_{L_{k-1}}v_1,\dots  ,A_{L_{k-1}}v_r\},$$
      $$\mathcal L'=\{L'_0,A_{L'_0}v'_1,\dots  ,A_{L'_0}v'_r,\dots ,
      L'_{k-1},A_{L'_{k-1}}v'_1,\dots  ,A_{L'_{k-1}}v'_r\}$$
      be two orthonormal bases constructed as above such that
      $L_0=X$, $L'_0=Y$, $g(L_\al,L_\al)=g(L'_\al,L'_\al)$ for
      $\al\in\{1,\dots  ,k-1\}$, and that
      $\{v_1=A_XY_1,\dots,v_r=A_XY_r\}$ and
      $\{v'_1=A_YY'_1,\dots,v'_r=A_YY'_r\}$ are orthonormal bases of the tangent
      bundle of the fibre $\pi^{-1}(\pi(p))$,
      where $Y_1,\dots,Y_r$ and $Y'_1,\dots,Y'_r$ are the horizontal liftings
      along $\pi^{-1}(\pi(p))$ of the vectors
      $\pi_*A_Xv_{1p},\dots,\pi_*A_Xv_{rp}$ and $\pi_*A_Yv'_{1p},\dots, \pi_*A_Yv'_{rp}$,
      respectively (as in Lemma \ref{lem:2}), for which
      $g(v_i,v_j)=g(v'_i,v'_j)$ for $i,j\in\{1,\cdots,r\}$.
      Let $\phi:T_pM\to T_pM$ be the linear map given by
      $\phi(L_\al)=L'_\al$, $\phi(v_j)=v'_j$, $\phi(A_{L_\al}v_j)=A_{L'_\al}v'_j$
      for every $\al\in\{0,\dots,k-1\}$ and
      $j\in\{1,\dots,r\}$.
      Since both $\mathcal L$, $\mathcal L'$ are orthonormal bases, we see that
      $\phi$ is a linear isometry.

      We shall apply Theorem \ref{t:4}.
      Thus we need to prove that $\phi(A_EF)=A_{\phi(E)}\phi(F)$ for every
      $E$, $F\in T_pM$.
      Indeed, we obtain for $\al, \beta\in\{0,\dots,k-1\}$
      and $j, l\in\{1,\dots,r\}$,
      \begin{eqnarray*}
            \phi(A_{L_\al}L_\bet)&=& \phi(0)=0=
            A_{L'_\al}L'_\bet=A_{\phi(L_\al)}\phi(L_\bet)\ ,\\
            g(v_j,A_{L_\al}A_{L_\beta}v_l)&=& -g(A_{L_\al}v_j,A_{L_\beta}v_l)
            =-g(L_\al,L_\beta)g(v_j,v_l)\\
            &=& -g(L'_\al,L'_\beta)g(v'_j,v'_l)=g(v'_j,A_{L'_\al}A_{L'_\beta}v'_l).
      \end{eqnarray*}
      Hence $\phi(A_{L_\al}A_{L_\beta}v_l)=A_{\phi(L_\al)}\phi(A_{L_\beta}v_l)$.
\ble\label{l:A}
      $A_{L_\al}v_j$ is a basic vector field along the
      fibre $\pi^{-1}(\pi(p))$ for every $1\leq j\leq r$ and $\al\geq 0$.
\ele
\begin{proof}[Proof of Lemma \rm{\ref{l:A}}]
      We have $g(A_Xv_j,Z)=g(A_XA_XY_j,Z)=-g(A_XY_j,A_XZ)$.
      For every  basic vector field $Z$ along the fibre $\pi^{-1}(\pi(p))$
      we know that $g(A_XY_j,A_XZ)$ is constant along the fibre $\pi^{-1}(\pi(p))$. Hence
      $A_Xv_j$ is a basic vector field along the fibre $\pi^{-1}(\pi(p))$.

      Now we assume $\al\geq 1$.
      Since $\dim(\ker A^*_X\cup\ker A^*_{L_\al})=\dim\ker A^*_X+\dim\ker A^*_{L_\al}-
      \dim(\ker A^*_X\cap\ker A^*_{L_\al})=(n-r)+(n-r)-(n-2r)=n$,
      it follows that $\ker A^*_X+\ker A^*_{L_\al}=\mathcal H$.
      Hence $A_{L_\al}v_j$ is a basic vector field along the
      fibre $\pi^{-1}(\pi(p))$ if and only if
      the following conditions are satisfied:
      $g(A_{L_\al}v_j,Z_1)$ is constant along $\pi^{-1}(\pi(p))$ for every
      $Z_1\in\ker A^*_X$, which is a  basic vector field along  $\pi^{-1}(\pi(p))$,
      and
      $g(A_{L_\al}v_j,Z_2)$ is constant along the fibre $\pi^{-1}(\pi(p))$ for every
      $Z_2\in\ker A^*_{L_\al}$, which is a basic vector field along $\pi^{-1}(\pi(p))$.
      If $Z_2\in\ker A^*_{L_\al}$, then $A^*_{L_\al}Z_2=0$ along $\pi^{-1}(\pi(p))$.
      So $g(A_{L_\al}v_j,Z_2)=-g(v_j,A_{L_\al}Z_2)=0$ along $\pi^{-1}(\pi(p))$.
      If $Z_1\in\ker A^*_X$, then $A^*_XZ_1=0$ along $\pi^{-1}(\pi(p))$.
      By O'Neill's equations, we get along the fibre $\pi^{-1}(\pi(p))$
      \begin{eqnarray*}
            R'(\pi_*X,\pi_*Y_j,\pi_*{L_\al},\pi_*Z_1)&=&
            R(X,Y_j,{L_\al},Z_1)+2g(A_XY_j,A_{L_\al}Z_1)\\
            &&-g(A_{Y_j}{L_\al},A_XZ_1)-g(A_{L_\al}X,A_{Y_j}Z_1)\\
            &=&-g(X,{L_\al})g(Y_j,Z_1)+g(X,Z_1)g(Y_j,{L_\al})\\
              &&+2g(v_j,A_{L_\al}Z_1),
      \end{eqnarray*}
      since $A_{L_\al}X=-A_X{L_\al}=0$ and $A_XZ_1=0$.
      Hence $g(v_j,A_{L_\al}Z_1)=-g(A_{L_\al}v_j,Z_1)$
      is constant along $\pi^{-1}(\pi(p))$
      for every
      $Z_1\in\ker A^*_X$, which is a basic vector field along $\pi^{-1}(\pi(p))$.

      We proved that $A_{L_\al}v_j$ is a basic vector field along $\pi^{-1}(\pi(p))$
      for every $\al\geq 0$ and $j\in\{1,\dots,r\}$.
\end{proof}
      We denote by $\hat\na$ the induced Levi-Civita connection on the fibre
      $\pi^{-1}(\pi(p))$.
\ble\label{l:B}
$A_{A_{L_\al}v_i}A_{L_\bet} v_j=g(L_\al,L_\bet)\hat\na_{v_i}v_j.$
\ele
\begin{proof}[Proof of Lemma \rm{\ref{l:B}}]
      By the relation \eqref{e:ecu} together with Lemma \ref{l:A}, we obtain for
      $i,j,l\in\{1,\dots,r\}$ and $\al,\beta\geq 0$ that
      \begin{eqnarray*}
            g(A_{A_{L_\al}v_i}{A_{L_\bet}v_j},v_l)&=& -g(A_{A_{L_\al}v_i}v_l,A_{L_\bet}v_j)\\
            &=&-R({L_\bet},v_l,A_{L_\al}v_i,v_j)+g(\na_{v_l}A_{L_\bet}A_{L_\al}v_i,v_j)\\
            &=&g({L_\bet},A_{L_\al}v_i)g(v_l,v_j)+v_lg(A_{L_\bet}A_{L_\al}v_i,v_j)\\
            &&-g(A_{L_\bet}A_{L_\al}v_i,\na_{v_l}v_j)\\
            &=&-v_lg(A_{L_\al}v_i,A_{L_\bet}v_j)+
               g(A_{L_\al}v_i,A_{L_\bet}v_t)g(\na_{v_l}v_j,v_t)\varepsilon_t\\
            &=&-g(L_\al,L_\bet)g(\hat\na_{v_l}v_j,v_i)\\
            &=&g(L_\al,L_\bet)g(\hat\na_{v_i}v_j,v_l).
      \end{eqnarray*}
      In the last equality we used the fact that $v_j=A_XY_j$
      is a Killing vector field along the fibre $\pi^{-1}(\pi(p))$
      (see \cite{bis} or \cite{bes}).
      Thus $$A_{A_{L_\al}v_i}{A_{L_\bet}v_j}=g(L_\al,L_\bet)\hat\na_{v_i}v_j.$$
\end{proof}
\ble\label{l:C}
The following assertions are true$:$
\begin{itemize}
\item[(a)] $r\not=2$.
\item[(b)] If $r=1$, then $A_{A_{L_\al}v_1}{A_{L_\bet}v_1}=0$.
\item[(c)] If $r=3$ and if we set $v_{3p}=(\hat\na_{v_1}v_2)(p)$, then $v_3=\hat\na_{v_1}v_2$
and
\begin{equation*}
     g(\hat\na_{v_i}v_j, v_k)=
        \left\{
          \begin{array}{cl}
              0 &  \mbox{if two of $i$, $j$, $k$ are equal,}\\
              \varepsilon \binom{1\ 2\ 3}{i\ j\ k}g(v_3,v_3)
                &  \mbox{if } \{i,j,k\}=\{1 , 2, 3\},
\end{array}
\right.
\end{equation*}
      where $\varepsilon \binom{1\ 2\ 3}{i\ j\ k}$ is the signature of the permutation $\binom{1\ 2\ 3}{i\ j\ k}$.
\end{itemize}
\ele
\begin{proof}[Proof of Lemma \rm{\ref{l:C}}]
      Since $v_1,\dots,v_r$ are Killing vector fields
      along $\pi^{-1}(\pi(p))$ and $g(v_i,v_i)\in\{-1,1\}$ for every $i$, we get
      $$g(\hat\na_{v_i}v_j,v_i)=g(\hat\na_{v_i}v_i,v_j)=g(\hat\na_{v_j}v_i,v_i)=0$$
      for every $i,j\in\{1,\dots,r\}$.

    (a)  The case $r=2$ is not possible. Indeed, if $r=2$, then the relation
      $g(\na_{v_1}v_2,v_1)=g(\na_{v_1}v_2,v_2)=0$ implies $\na_{v_1}v_2=0$.
      On the other hand,
      $$g(\na_{v_1}v_2,\na_{v_1}v_2)=-g(\hat\na_{v_1}\hat\na_{v_2}v_2,v_1)+
      \hat R(v_1,v_2,v_1,v_2)=-g(v_1,v_1)g(v_2,v_2)\in\{-1,1\},$$
      since $\hat\na_{v_2}v_2=g(X,X)^{-1}A_{A_Xv_2}A_Xv_2=0$
      and each fibre has constant curvature $-1$.
      So we get a contradiction.

    (b)  If  $r=1$, then $A_{A_{L_\al}v_1}A_{L_\bet}v_1=0$
      for every $\al$ and $\bet$, because
      $0=A_{A_Xv_1}A_Xv_1=g(X,X)\na_{v_1}v_1$ implies $\na_{v_1}v_1=0$.

    (c)  In the case $r=3$ we shall prove
      $g(\hat\na_{v_1}v_2,v_3)$ is constant along the fibre $\pi^{-1}(\pi(p))$.
      Since O'Neill's integrability tensor $A$ is skew-symmetric, it follows that
      $\hat\na_{v_i}v_j=-\hat\na_{v_j}v_i$.
      Then $\hat\na_{v_i}v_j=(1/2)[v_i,v_j]$ is a Killing vector field along
        $\pi^{-1}(\pi(p))$.
      We then obtain
      \begin{eqnarray*}
            v_1g(\hat\na_{v_1}v_2,v_3)&=&
            g(\hat\na_{v_1}\hat\na_{v_1}v_2,v_3)+g(\hat\na_{v_1}v_2,\hat\na_{v_1}v_3)\\
            &=&-g(\hat\na_{v_3}\hat\na_{v_1}v_2,v_1)+g(\hat\na_{v_1}v_2,\hat\na_{v_1}v_3)\\
            &=&-v_3g(\hat\na_{v_1}v_2,v_1)+g(\hat\na_{v_1}v_2,\hat\na_{v_1}v_3+\hat\na_{v_3}v_1)=0.
            \end{eqnarray*}
      Analogously, we get $v_2g(\hat\na_{v_1}v_2,v_3)=-v_2g(\hat\na_{v_2}v_1,v_3)=0$. We also obtain
      $$v_3g(\hat\na_{v_1}v_2,v_3)=
      g(\hat\na_{v_3}\hat\na_{v_1}v_2,v_3)+g(\hat\na_{v_1}v_2,\hat\na_{v_3}v_3)=0,$$
      since $\hat\na_{v_3}v_3=0$ and $\hat\na_{v_1}v_2$
      is a Killing vector field along  $\pi^{-1}(\pi(p))$.
      It is easy to see that
\begin{eqnarray*}
        g(\hat\na_{v_1}v_2,v_3) &=& -g(\hat\na_{v_2}v_1,v_3)=g(\hat\na_{v_2}v_3,v_1)\\
        &=&-g(\hat\na_{v_3}v_2,v_1)=g(\hat\na_{v_3}v_1,v_2)=-g(\hat\na_{v_1}v_3,v_2).
\end{eqnarray*}
      Thus $g(\hat\na_{v_i}v_j,v_l)$ is constant along the fibre $\pi^{-1}(\pi(p))$ for each
      $i,j,l\in\{1,2,3\}$.
      Therefore
      $$g(A_XA_{A_Xv_i}A_Xv_j,A_Xv_l)=-g(X,X)g(A_{A_Xv_i}A_Xv_j,v_l)
       =-g(X,X)^2g(\hat\na_{v_i}v_j,v_l)$$
      is constant along $\pi^{-1}(\pi(p))$.
      Also, we compute for $\al\geq 1$
    $$g(A_XA_{A_Xv_i}A_Xv_j, A_{L_\al}v_l)=-g(A_{A_Xv_i}A_Xv_j, A_XA_{L_\al}v_l)=0,$$
    $$g(A_XA_{A_Xv_i}A_Xv_j,L_\al)=-g(A_{A_Xv_i}A_Xv_j,A_XL_\al)=0.$$
      Hence $A_XA_{A_Xv_i}A_Xv_j=g(X,X)A_X\hat\na_{v_i}v_j$
      is a basic vector field for each $i,j\in\{1,\dots,r\}$.

      We choose $v_{3p}=(\hat\na_{v_1}v_2)(p)$.
      Since $A_X\hat\na_{v_1}v_2$ is a basic vector
      field along  $\pi^{-1}(\pi(p))$, we get the horizontal lifting along
      $\pi^{-1}(\pi(p))$ of
      $\pi_*(A_X\hat\na_{v_1}v_2(p))=\pi_*A_Xv_{3p}$ is
      $g(X,X)^{-1}A_X\hat\na_{v_1}v_2$.
      On the other hand, $Y_3$ is, by definition, the horizontal lifting of $g(X,X)^{-1}\pi_*A_Xv_{3p}$ along
      $\pi^{-1}(\pi(p))$.  It follows that
       $Y_3=g(X,X)^{-1}A_X\hat\na_{v_1}v_2$
      along $\pi^{-1}(\pi(p))$.
      Thus
      $$v_3=A_XY_3=g(X,X)^{-1}A_XA_X\hat\na_{v_1}v_2=\hat\na_{v_1}v_2$$
      along the fibre $\pi^{-1}(\pi(p))$.
\end{proof}

      For $r=3$, we choose $v'_{3p}=(\hat\na_{v'_1}v'_2)(p)$.
      If we repeat the argument above for the
      basis $\{v'_1,v'_2,v'_3\}$, by Lemma \ref{l:C}, we get
      $v'_3=\hat\na_{v'_1}v'_2$ along the fibre $\pi^{-1}(\pi(p))$.
      It follows that
      $g(\hat\na_{v_i}v_j,v_l)= g(\hat\na_{v'_i}v'_j,v'_l)$
      for each $i,j,l\in\{1,2,3\}$.

      Returning to the computation of
      $g(A_{A_{L_\al}v_i}A_{L_\bet} v_j,v_l)$,
      in both cases $r=1$ and $r=3$, we get
      for every $\al,\bet\geq 0$ and $i,j,k\in\{1,\dots,r\}$
      \begin{eqnarray*}
            g(A_{A_{L_\al}v_i}A_{L_\bet} v_j,v_l)&=& g(L_\al,L_\bet)g(\hat\na_{v_i}v_j,v_l)\\
      &=& g(L'_\al,L'_\bet)g(\hat\na_{v'_i}v'_j,v'_l)=g(A_{A_{L'_\al}v'_i}A_{L'_\bet} v'_j,v'_l).
      \end{eqnarray*}
      Hence
      $\phi(A_{A_{L_\al}v_i}A_{L_\bet} v_j)=A_{\phi(A_{L_\al}v_i)}\phi(A_{L_\bet} v_j)$ and
      $\phi(A_{A_{L_\al}v_i}v_j)=A_{\phi(A_{L_\al}v_i)}\phi(v_j).$

      By Corollary 2.3.14 in \cite{wol} we see that $\phi:T_pM\to T_pM$
      extends to an isometry on $M$, denoted by $f:M\to M$, such that
      $f(p)=p$  and $f_{*p}=\phi$. Hence $f_{*p}X=Y$ and
      $f_*(\mathcal H_p)=\mathcal H_p$.
      Since $f_{*}A_EF=A_{f_*E}f_*F$ for every $E, F\in T_pM$, we see, by Theorem
      \ref{t:4}, that there is an isometry  $\Tilde{f}:B\to B$ such that
      $\Tilde{f}\circ\pi=\pi\circ f$.
      Thus $\Tilde f_*X'=\Tilde f_*\pi_*X=\pi_*f_*X=\pi_*Y=Y'$ and
      $\Tilde f(x)=\Tilde f(\pi(p))=\pi(f(p))=\pi(p)=x$.

      Therefore $B$ is an isotropic semi-Riemannian manifold.
      This completes the proof of Theorem \ref{p:3}.
\end{proof}
      If the metric on the base space is negative definite,
       the following lemma follows from Theorem \ref{p:3}.
\ble\label{l:new}
    If $\pi:M\to B$ is a semi-Riemannian submersion
    with  connected totally geodesic fibres
    from an $(n+r)$-dimensional semi-Riemannian manifold $M$
    of index $r'+n$ and of  constant  negative curvature
    onto an $n$-dimensional  semi-Riemannian manifold $B$ of index $n$, then
    $r'=r$.
\ele
\begin{proof}
      By Theorem \ref{p:3}, we have $n=q_1(r'+1)+q_2(r-r')=(q_1+q_2)(r+1)$ for
      some nonnegative integers  $q_1$ and  $q_2$. Hence
      $0=q_1(r-r')+q_2(r'+1)$. Since the right hand side
      is the sum of two non-negative numbers, it follows that
      $q_1(r-r')=0$ and $q_2(r'+1)=0$. Therefore $q_2=0$. This implies $r'=r$.
\end{proof}
\begin{remark}
    Changing simultaneously the signs of metrics
    on the total space and on the base space,
    any semi-Riemannian submersion, under the assumptions of Lemma \ref{l:new},
    becomes a Riemannian submersion
    with totally geodesic fibres from a sphere onto a Riemannian manifold.
    This case was completely classified by Escobales (see \cite{esc}) and
    Ranjan (see \cite{ran}).
\end{remark}
\bpr\label{p:5}
    Let $\pi:M\to B$ be a semi-Riemannian submersion
    with connected totally geodesic fibres from a complete simply connected
    semi-Riemannian manifold $M$ onto a semi-Riemannian manifold $B$. Then
    $B$ is simply connected and complete.
\epr
\begin{proof}
      If $M$ is geodesically complete, then so is $B$ (see
      \cite{bes} or \cite{bad}).
      Since $M$ is a complete semi-Riemannian manifold  and the fibres are
      totally geodesic, any fibre is also geodesically complete. By a theorem in
      \cite{rec}, it follows that the horizontal distribution $\mathcal H$  is an Ehresmann
      connection. Therefore, by \cite{ehr}, we see that $\pi$ is a fibre bundle. So we obtain
      an exact homotopy sequence:
      $$\dots\to\pi_2(M)\to\pi_2(B)\to\pi_1(fibre)\to\pi_1(M)\to\pi_1(B)\to 0.$$
      Thus $\pi_1(B)=0$.
\end{proof}

By Theorem 12.3.2 in \cite{wol}, we know that any
connected, simply connected isotropic semi-Riemannian
manifold is isometric to one of the following semi-Riemannian
manifolds:
\begin{itemize}
      \item[(i)] $\mathbb R^m_t$ or the universal semi-Riemannian
                 covering of the  pseudo-hyperbolic space  $H^m_t(c)$ with constant  sectional
      curvature $c<0$, or of the  pseudo-sphere $S^m_t(c)$ with constant  sectional
      curvature $c>0$.
\item[(ii)] The complex pseudo-hyperbolic
      space  $\mathbb CH^m_t(c)$ with constant holomorphic sectional
      curvature $c<0$, or
       the complex pseudo-projective space $\mathbb CP^m_t(c)$ with constant holomorphic sectional
      curvature $c>0$.
\item[(iii)] The quaternionic pseudo-hyperbolic space $\mathbb HH^m_t(c)$
      with constant quaternionic sectional
      curvature $c<0$,
      or the quaternionic pseudo-projective space  $\mathbb HP^m_t(c)$
      with constant quaternionic sectional
      curvature $c>0$.
\item[(iv)] The Cayley pseudo-hyperbolic plane  $\mathbb CaH^2_t(c)$ with Cayley sectional
      curvature $c<0$,
      or the Cayley pseudo-projective plane  $\mathbb CaP^2_t(c)$  with Cayley sectional
      curvature $c>0$.
\end{itemize}

\ble
  {\rm (a)} If $B$ is a semi-Riemannian manifold isometric to one
       of the semi-Riemannian manifolds
       $\mathbb CP^m_t(c)$,  $\mathbb HP^m_t(c)$, $\mathbb CaP^2_t(c)$ $(c>0)$,
       then the curvature tensor satisfies the inequality
\begin{equation}\label{e:*}
       R'(X',Y',X',Y')\geq\f c4 ( g'(X',X')g'(Y',Y')-g'(X',Y')^2)
\end{equation}
       for each tangent vectors  $X'$, $Y'$ of $B$.

  {\rm (b)} If $B$ is a semi-Riemannian manifold isometric to one
       of the semi-Riemannian manifolds
       $\mathbb CH^m_t(c)$,  $\mathbb HH^m_t(c)$, $\mathbb CaH^2_t(c)$ $(c<0)$,
       then the curvature tensor satisfies the inequality
\begin{equation}\label{e:**}
       R'(X',Y',X',Y')\leq\f c4 ( g'(X',X')g'(Y',Y')-g'(X',Y')^2)
\end{equation}
       for each tangent vectors $X'$, $Y'$  of $B$.
\ele
\begin{proof}
    For each  tangent vectors $X'$, $Y'$ of  $B$, we have the following formulas
    for the
    curvature tensors:
\begin{itemize}
\item[(i)] If $B\in\{\mathbb CP^m_t(c),\mathbb CH^m_t(c)\}$ and
    $I_0$ is the natural complex structure on $B$, then
\begin{equation}\label{e:i}
R'(X',Y',X',Y')=\f c4 (g'(X',X')g'(Y',Y')-g'(X',Y')^2+3g'(X',I_0Y')^2).
\end{equation}
\item[(ii)] If $B\in\{\mathbb HP^m_t(c),\mathbb HH^m_t(c)\}$ and
    $I_0,J_0,K_0$ are local almost complex structures
    which give rise to the quaternionic structure on $B$, then
\begin{eqnarray}\label{e:ii}
    R'(X',Y',X',Y')&=&(c/4)(g'(X',X')g'(Y',Y')-g'(X',Y')^2 \\
   &&+3g'(X',I_0Y')^2+3g'(X',J_0Y')^2+3g'(X',K_0Y')^2).\no
\end{eqnarray}
\end{itemize}
By these explicit formulas for curvature tensors, in all cases we obtain
the inequalities \eqref{e:*} and \eqref{e:**}.
\end{proof}

First, we shall discuss the case of a base space with nonconstant curvature.
\ble\label{l:6}
    If $\pi:H^{n+r}_{s+r'}\to B^n_s$ is a semi-Riemannian submersion
        with  connected totally geodesic fibres from an
    $(n+r)$-dimensional pseudo-hyperbolic space $H^{n+r}_{s+r'}$  of index $s+r'>1$
    onto an $n$-dimensional isotropic  semi-Riemannian
    manifold  $B^n_s$ of index $s$ with nonconstant
    curvature, then the induced metrics on the fibres are negative definite
    and $B$ is isometric to one of the following semi-Riemannian manifolds$:$
    \begin{itemize}
        \item[{\rm (i)}] $\mathbb CH^m_t$, $m>1$,
        \item[{\rm (ii)}] $\mathbb HH^m_t$, $m>1$,
        \item[{\rm (iii)}] $\mathbb CaH^2_t$.
    \end{itemize}
\ele
\begin{proof}
      Since $\dim\mathcal H=k(\dim\mathcal V+1)$ for some  positive integer $k$,
      we get $\dim \mathcal H\geq \dim \mathcal V+1$.
      Let $X$ be a horizontal vector field along a fibre $\pi^{-1}(\pi(p))$ such that
      $g(X,X)\neq 0$ and $X$ is the horizontal lifting of some tangent vector of $B$.

      First, we shall prove that
      $$\dim \mathcal H> \dim \mathcal V+1.$$
      Suppose  that
      $\dim \mathcal H= \dim \mathcal V+1$.
      Then
      $A_X:\mathcal V\to X^\perp=\{Y\in\mathcal H\ |\ g(X,Y)=0\}$
      is bijective. For every $Y\in X^\perp$ we get $Y=A_XV$ for some
      vertical vector $V$. It follows that
\begin{eqnarray*}
      g(A_XY,A_XY)&=& g(A_XA_XV,A_XA_XV)=g(X,X)^2g(V,V), \\
      g(Y,Y)&=& g(A_XV,A_XV)=-g(X,X)g(V,V).
\end{eqnarray*}
      Thus
      $g(A_XY,A_XY)=-g(X,X)g(Y,Y)$ for every $Y\in X^\perp$.
      By O'Neill's equations, we have
\begin{eqnarray*}
      R'(\pi_*X,\pi_*Y,\pi_*X,\pi_*Y)&=& -g(X,X)g(Y,Y)+g(X,Y)^2+3g(A_XY,A_XY)\\
      &=& -4(g(X,X)g(Y,Y)-g(X,Y)^2)
\end{eqnarray*}
      for every horizontal vector field  $Y$ along
      $\pi^{-1}(\pi(p))$.
      Hence $B$ has constant curvature is a contradiction.

      We established that $\dim\mathcal H> \dim\mathcal V+1.$
      So we can find  a horizontal vector field $Z$
      along the fibre $\pi^{-1}(\pi(p))$
      such that $Z\in\ker A^*_X$, $g(X,Z)=0$, $g(Z,Z)\neq 0$ and
      $Z$ is the horizontal lifting of some $Z'\in T_{\pi(p)}B$.
      We then have
\begin{eqnarray*}
      R'(\pi_*X,\pi_*Z,\pi_*X,\pi_*Z)&=& -g(X,X)g(Z,Z)+g(X,Z)^2+3g(A_XZ,A_XZ)\\
      &=& -g(X,X)g(Z,Z).
\end{eqnarray*}
      Since $B$ is a simply connected isotropic semi-Riemannian manifold
      with nonconstant curvature,
      we see that  $B$ is isometric to one of the following semi-Riemannian manifolds:
\begin{itemize}
      \item[(a) ] $\mathbb CP^m_t(c)$, $\mathbb HP^m_t(c)$, $\mathbb CaP^2_t(c)$, or
      \item[(b) ] $\mathbb CH^m_t(c)$, $\mathbb HH^m_t(c)$, $\mathbb CaH^2_t(c)$.
\end{itemize}
We shall prove that only the case (b) is possible.

First, we suppose that $B$ is isometric to one of the following semi-Riemannian manifolds:
\begin{center}
       $\mathbb CP^m_t(c)$, $\mathbb HP^m_t(c)$, $\mathbb CaP^2_t(c)$ $(c>0)$.
\end{center}
      By the inequality \eqref{e:*}, we get
\begin{eqnarray*}
      R'(\pi_*X,\pi_*A_XV,\pi_*X,\pi_*A_XV)&=&-4g(X,X)g(A_XV,A_XV)\\
      &=&4g(X,X)^2g(V,V)\geq -(c/4) g(X,X)^2g(V,V).
\end{eqnarray*}
      Therefore
\begin{equation}\label{e:t1}
      g(V,V)\geq 0
\end{equation}
      for every  vertical vector $V$.
      Since $X$ and $Z$ are basic vector fields along
      $\pi^{-1}(\pi(p))$ with  $g(X,Z)=0$ and $A_XZ=0$ along
      $\pi^{-1}(\pi(p))$, it follows from the relation \eqref{e:ecu} that $A_ZV\in\ker A^*_X$.
      On the other hand, by the inequality \eqref{e:*}, we get
\begin{eqnarray*}
      R'(\pi_*X,\pi_*Z,\pi_*X,\pi_*Z)&=& -g(X,X)g(Z,Z)\geq (c/4)  g(X,X)g(Z,Z),\\
      R'(\pi_*X,\pi_*A_ZV,\pi_*X,\pi_*A_ZV)&=& -g(X,X)g(A_ZV,A_ZV)\geq (c/4) g(X,X)g(A_ZV,A_ZV).
\end{eqnarray*}
      Hence $g(X,X)g(Z,Z)\leq 0$ and $g(X,X)g(A_ZV,A_ZV)\leq 0$.
      Thus
            $$0\leq g(Z,Z)g(A_ZV,A_ZV)=-g(Z,Z)^2g(V,V).$$
      So for any  vertical vector $V$ we get
\begin{equation}\label{e:t2}
      g(V,V)\leq 0.
\end{equation}
      Since the induced metrics on fibres are nondegenerate, it is not possible
      to have both \eqref{e:t1} and \eqref{e:t2}.
      So we obtain the required contradiction.
      It follows that $B$ is isometric to one of the following semi-Riemannian
      manifolds:
\begin{center}
       $\mathbb CH^m_t(c)$, $\mathbb HH^m_t(c)$, $\mathbb CaH^2_t(c)\ $ $(c<0)$.
\end{center}

      We shall now prove that $c=-4$. Suppose $(c/4)  +1\not=0$.
      By the inequality \eqref{e:**}, we get
\begin{equation}\label{e:***}
      R'(\pi_*X,\pi_*Z,\pi_*X,\pi_*Z)= -g(X,X)g(Z,Z)\leq (c/4) g(X,X)g(Z,Z),
\end{equation}
      $$R'(\pi_*X,\pi_*A_ZV,\pi_*X,\pi_*A_ZV)= -g(X,X)g(A_ZV,A_ZV)\leq (c/4)g(X,X)g(A_ZV,A_ZV).$$
      Hence
\begin{equation}\label{e:****}
      ((c/4)+1)^2g(X,X)^2g(Z,Z)g(A_ZV,A_ZV)\geq 0,
\end{equation}
      from which follows that $0\leq g(Z,Z)g(A_ZV,A_ZV)=-g(Z,Z)^2g(V,V)$.
     Therefore $g(V,V)\leq 0$ for every
      vertical vector field $V$. In particular, we have $g(A_XY,A_XY)\leq 0$, which implies
\begin{equation}\label{e:9}
      R'(\pi_*X,\pi_*Y,\pi_*X,\pi_*Y)\leq g(X,X)g(Y,Y)-g(X,Y)^2
\end{equation}
      for every  horizontal vectors $X$ and $Y$.
      We have the following cases:

\vspace*{5pt}
{\it Case} (a)\ $0<\mathrm{index\ }B<\dim B$.
      We can choose vector fields  $X'$, $Y'$  on $B$ such that
     $g'(X',X')g'(Y',Y')<0$ and
      that one of the following conditions is satisfied:
\begin{itemize}
\item[(i)] $Y'\in\{X',I_0X'\}^\perp$  if $B=\mathbb CH^m_s(c)$, where
          $I_0$ is the natural complex structure on $\mathbb CH^m_s(c)$,
\item[(ii)] $Y'\in\{X',I_0X',J_0X',K_0X'\}^\perp$  if $B=\mathbb HH^m_s(c)$, where
          $\{I_0, J_0, K_0\}$ are local almost complex structures
          which give rise to the quaternionic structure on
          $\mathbb HH^m_s(c)$, or
\end{itemize}
      Let $X$, $Y$ be the horizontal liftings of $X'$, $Y'$. The inequality \eqref{e:9}
      then implies
      $$\f c4g(X,X)g(Y,Y)\leq -g(X,X)g(Y,Y).$$
      Hence
         $((c/4)+1)g(X,X)g(Y,Y)\leq 0$. Therefore $(c/4)+1>0$.
      On the other hand, we can choose horizontal vector fields  $X$, $Z$ such that
      $g(X,Z)=0$, $Z\in\ker A^*_X$ and $g(X,X)g(Z,Z)<0$, because
            $0<\mathrm{index\ }B<\dim B$. Then the inequality \eqref{e:***} becomes
      $(c/4)+1<0$. So we get a contradiction.

\vspace*{5pt}
{\it Case} (b)\  $\mathrm{index\ }B\in\{0,\dim B\}$.
      Similarly, we can choose vector fields $X'$, $Y'$ on $B$ such that
      $g'(X',Y')=0$ and $R'(X',Y',X',Y')=(c/4)g'(X',X')g'(Y',Y')$.
      The inequality \eqref{e:9} then implies
      $((c/4)+1)g'(X',X')g'(Y',Y')\leq 0$. By the hypothesis
      of Case (b), we get $(c/4)+1\leq 0$. On the other hand, the inequality
      \eqref{e:***} becomes $(c/4)+1>0$. So we get a contradiction.

\vspace*{7pt}
      We have proved $c=-4$. The inequality \eqref{e:**} then becomes
\begin{equation}
      R'(X',Y',X',Y')\leq -g'(X',X')g'(Y',Y')+g'(X',Y')^2
\end{equation}
      for tangent vector fields $X'$, $Y'$ on $B$.
      Then we have
      $$R'(\pi_*X,\pi_*A_XV,\pi_*X,\pi_*A_XV)=-4g(X,X)g(A_XV,A_XV)\leq -g(X,X)g(A_XV,A_XV)$$
      for a vertical vector field  $V$ and for a  horizontal vector field $X$.
      Hence
      $$0\leq g(X,X)g(A_XV,A_XV)=-g(X,X)^2g(V,V).$$
      Therefore
      the induced metrics on fibres are negative definite.
\end{proof}

By Lemma \ref{l:6}, we deduce the following proposition.
\bpr\label{p:8}
    If $\pi:H^{n+r}_{s+r'}\to B^n_s$ is a semi-Riemannian submersion with
       connected  totally geodesic fibres from an
    $(n+r)$-dimensional pseudo-hyperbolic space $H^{n+r}_{s+r'}$ of index $s+r'$
    onto an $n$-dimensional isotropic  semi-Riemannian
    manifold  $B^n_s$ of index $s$ with nonconstant
    curvature, and if the fibres are negatively definite
    then one of the following holds$:$
    \begin{itemize}
        \item[{\rm (1)}] $n=2m>2$, $s=2t$, $r=r'=1$
                 for some  non-negative integers $m$, $t$,
        and $B^n_s$ is isometric to $\mathbb CH^m_t$.
        \item[{\rm (2)}] $n=4m>4$, $s=4t$, $r=r'=3$
                 for some non-negative integers $m$, $t$,
        and $B^n_s$ is isometric to $\mathbb HH^m_t$.
        \item[{\rm (3)}] $n=16$, $s\in\{0, 8, 16\}$, $r=r'=7$, and
        $B^n_s$ is isometric to $\mathbb CaH^2_{s/8}$.
    \end{itemize}
\epr
\begin{proof}
First, we shall discuss the case $s+r'>1$.
By Lemma \ref{l:6}, $B$ is isometric
to one of the semi-Riemannian manifolds
$\mathbb CH^m_t$, $\mathbb HH^m_t$, $\mathbb CaH^2_t$ for some $m>1$.

Let $x\in B$ and let $X'\in T_xB$ such that $g'(X',X')\neq 0$,
and let $\mathcal F_{X'}$
be the subspace in $T_xB$ given by
$$\mathcal F_{X'}=\{Y'\in T_xB\ |\ R'(X',Y')X'=-g'(X',Y')X'+g'(X',X')Y'\}.$$
Let $p\in\pi^{-1}(x)$ and let $X$ be the horizontal
lifting vector at $p$ of $X'$.
By O'Neill's equations, we have
$R'(\pi_*X,\pi_*Y,\pi_*X,\pi_*Z)=R(X,Y,X,Z)+3g(A^*_XY,A^*_XZ)$
for  horizontal vectors $Y$, $Z$.
Since $A^*_X:\mathcal H_p\to\mathcal V_p$ is surjective
and since the induced metrics on fibres are nondegenerate, we get
$Y\in\ker A^*_X$ if and only if $\pi_*Y\in\mathcal F_{X'}$.
Thus
$$\dim\ker A^*_X=\dim\mathcal H-\dim\mathcal V=\dim\mathcal F_{X'}.$$
We have the following possibilities:
\begin{itemize}
\item[(1) ] $B^n_s$ is isometric to $\mathbb CH^m_t$. So $n=2m$, $s=2t$.
         From the geometry of the complex
    pseudo-hyperbolic space (see relation \eqref{e:i}), we get
    $\dim \mathcal F_{X'}=\dim \mathcal H-1$.
    It follows that $r=r'=\dim \mathcal V=1$.
\item[(2) ] $B^n_s$ is isometric to $\mathbb HH^m_t$. So $n=4m$, $s=4t$.
         From the geometry of the quaternionic
    pseudo-hyperbolic space (see relation \eqref{e:ii}), we get
    $\dim \mathcal F_{X'}=\dim \mathcal H-3$.
    It follows that $r=r'=\dim \mathcal V=3$.
\item[(3) ] $B^n_s$ is isometric to  the Cayley pseudo-hyperbolic plane
     $\mathbb CaH^2_t$. So $n=16$, ${s\in\{0, 8, 16\}}$.
         From the geometry of the Cayley pseudo-hyperbolic plane, we obtain
     $\dim \mathcal F_{X'}=\dim \mathcal H-7$.
     Hence $r=r'=\dim \mathcal V=7$.
\end{itemize}
\vspace*{0.1cm}

Now, we discuss the remaining case $s+r'=1$.
From  $s+r'=1$, we have either
\begin{itemize}
\item[(i) ] $s=0$, $r'=1$,\ or
\item[(ii) ] $s=1$, $r'=0$.
\end{itemize}

If $s=0$, $r'=1$, then $\pi:H^{n+r}_1\to B^n$ is a semi-Riemannian submersion
with totally geodesic fibres from an anti-de Sitter space onto a Riemannian manifold.
In this case, investigated by Magid in \cite{mag}, it follows that $B$ is isometric to the
complex hyperbolic space  $\mathbb CH^m$ and $r=r'=1$.

For $s=1$, $r'=0$, we get, by Theorem \ref{p:3},
$1=q_1+q_2r\geq q_1+q_2$ with $q_1+q_2=k=n/(r+1)$.
Thus $q_1+q_2=1$. It follows that $n=r+1$. Hence $A_X:\mathcal V\to X^\perp$
is bijective. Since $R'(\pi_*X,\pi_*A_XV,\pi_*X,\pi_*A_XV)=-4g(X,X)g(A_XV,A_XV)$,
we see that $B$ has constant curvature $-4$,
which contradicts our assumption of
nonconstant curvature of the base space.
\end{proof}

We shall now discuss the case where the base space is of constant
curvature.

\bpr\label{p:9}
        If $\pi:H^{n+r}_{s+r'}\to B^n_s$ is a semi-Riemannian submersion
        with connected totally geodesic fibres from
        an $(n+r)$-dimensional pseudo-hyperbolic space of index $s+r'$ onto
        an $n$-dimensional semi-Riemannian manifold of index $s$ with constant curvature,
        and if the fibres are negatively definite,
        then one of the following holds$:$
        \begin{itemize}
                \item[(1)]  $n=s=2^t$, $r=r'=n-1$,
        $B$ is isometric to $H^{2^t}_{2^t}(-4)$ and
        $t\in\{1,2,3\}$.
                \item[(2)]  $n=2^t$, $s=0$, $r=r'=n-1$,
           $B$ is isometric to $H^{2^t}(-4)$ and
        $t\in\{1,2,3\}$.
        \end{itemize}
\epr
\begin{proof}
Since $B$ has constant curvature, the curvature of $B$ is $-4$ and
$n=r+1$. By Theorem \ref{p:3}, $s=q_1(r'+1)+q_2(r-r')=q_1(r+1)$ and
$q_1+q_2=n/(r+1)=1$. Then either $q_1=0$ or $q_1=1$. If $q_1=0$,
then $s=0$. If $q_1=1$ then $s=r+1=n$. Summarizing, we have
$\mathrm{index}(B)\in\{0,\dim B\}$.

If $\mathrm{index}(B)=\dim B$, then, by Lemma \ref{l:new}, we obtain $r=r'$.
Hence, by \cite{ran}, we have (1).

If $\mathrm{index}(B)=0$, then, by \cite{bad}, we have (2).\\
The idea of the proof in \cite{ran} and \cite{bad} is to see that the
tangent bundle of any fibre is trivial and that fibres are
diffeomorphic to spheres, and then to apply a well-known result of
Adams which claims that the spheres of dimensions 1, 3 and 7 are the
only spheres with trivial tangent bundle.
\end{proof}

The next theorems solve the equivalence problem
of  semi-Riemannian submersions from
real and complex pseudo-hyperbolic spaces.
\bth\label{t:eq}
    If $\pi_1,\pi_2:H^{n+r}_{s+r'}\to B^n_s$ are
    two semi-Riemannian submersions with connected totally
    geodesic fibres from a pseudo-hyperbolic space of index
    $s+r'>1$, if the fibres are negative definite, and if
    the dimension of the fibres is $r\in\{1,3\}$,
    then $\pi_1$ and $\pi_2$ are equivalent.
\end{theorem}
\begin{proof}
Let $p,q\in H^{n+r}_{s+r'}$. Let
$$\mathcal L=\{L_0,A_{1L_0}v_1,\dots  ,A_{1L_0}v_r,\dots ,
L_{k-1},A_{1L_{k-1}}v_1,\dots  ,A_{1L_{k-1}}v_r\},$$
$$\mathcal L'=\{L'_0,A_{2L'_0}v'_1,\dots  ,A_{2L'_0}v'_r,\dots ,
L'_{k-1},A_{2L'_{k-1}}v'_1,\dots  ,A_{2L'_{k-1}}v'_r\}$$
be two orthonormal bases of $\mathcal H_{1}$ along $\pi_1^{-1}(\pi_1(p))$
 and of $\mathcal H_{2}$ along $\pi_2^{-1}(\pi_2(q))$ constructed as
in the proof of Theorem \ref{p:3} such that
$g_p(L_\al,L_\bet)=g_q(L'_\al,L'_\bet)=\varepsilon_\al\de_{\al\bet}$ for
$\al, \bet\in\{0,\dots  ,k-1\}$,
$g_p(v_i,v_j)=g_q(v'_i,v'_j)=\ei\de_{ij}$ for  $i,j\in\{1,\dots,r\}$ and
for $r=3$, $v_{3p}=(\hat\na_{v_1}v_2)(p)$ and $v'_{3q}=(\hat\na_{v'_1}v'_2)(q)$.

Let $\phi:T_pH^{n+r}_{s+r'}\to T_qH^{n+r}_{s+r'}$ be the linear map given by
$\phi(L_\al)=L'_\al$,
$\phi(A_{1L_\al}v_i)=A_{2L'_\al}v'_i$, $\phi(v_i)=v'_i$
for every $\al$ and $i$.
In a manner similar to the proof of Theorem \ref{p:3},
we obtain $\phi(A_{1E}F)=A_{2\phi(E)}\phi(F)$
for every $E$, $F\in T_pH^{n+r}_{s+r'}$.
By Corollary 2.3.14 in \cite{wol}, $\phi$ extends to an isometry on $H^{n+r}_{s+r'}$,
denoted by
$f:H^{n+r}_{s+r'}\to H^{n+r}_{s+r'}$, satisfying $f(p)=q$ and $f_{*p}=\phi$.
From Theorem \ref{t:4} it follows that
$f$ induces an isometry $\Tilde f$ on $B$,
such that $\tilde f\circ\pi=\pi\circ f$.
Hence $\pi_1$ and $\pi_2$ are equivalent.
\end{proof}
\bth\label{t:eq2}
    If $\pi_1,\pi_2:\mathbb CH^{2n+1}_{2s+1}\to \mathbb HH^n_s$ are
    two semi-Riemannian submersions with connected complex totally
    geodesic fibres from a complex pseudo-hyperbolic space, and if
    the fibres are negative definite,
    then $\pi_1$ and $\pi_2$ are equivalent.
\end{theorem}
\begin{proof}
Let
$\theta:H^{4n+3}_{4s+3}\to\mathbb CH^{2n+1}_{2s+1}$
be the canonical semi-Riemannian
submersion. By Theorem 2.5 in \cite{esco}, we see that
$\tilde\pi_1=\pi_1\circ\theta:H^{4n+3}_{4s+3}\to\mathbb HH^n_s$ and
$\tilde\pi_2=\pi_2\circ\theta:H^{4n+3}_{4s+3}\to\mathbb HH^n_s$
are semi-Riemannian submersions with totally geodesic fibres. We denote by
$\tilde A_1$, $\tilde A_2$, $A_1$, $A_2$, $A$  O'Neill's integrability tensors of
$\tilde \pi_1$, $\tilde \pi_2$, $\pi_1$, $\pi_2$, $\theta$, respectively.
In order to reduce the proof of the equivalence theorem of semi-Riemannian submersions
from a complex pseudo-hyperbolic space to that
from a pseudo-hyperbolic space, we need to establish relations among the integrability
tensors $\tilde A_1$, $A_1$, $A$.

First, we prove that $\theta_*\Tilde A_{1X}Y=A_{1\theta_*X}\theta_*Y$
for $\Tilde\pi_1$-basic vector fields $X$ and $Y$.
Let $p\in H^{4n+3}_{4s+3}$. Let $w'_1$, $w'_2$ be two orthonormal
$\pi_1$-vertical vectors in $T_{\theta(p)}\mathbb
CH^{2n+1}_{2s+1}$ and let $w_1$, $w_2$ be the $\theta$-horizontal
liftings at $p$ of $w'_1$, $w'_2$, respectively. Let $w_3$ be a
unit $\theta$-vertical vector in $T_pH^{4n+3}_{4s+3}$. Then
$\{w_1,w_2,w_3\}$ gives an orthonormal basis of $\tilde{\mathcal
V}_{1p}$.
Since the induced metrics on the fibres of $\tilde\pi_1$ are negative definite, we have
$$\Tilde A_{1X}Y=-g(\na_XY,w_1)w_1-g(\na_XY,w_2)w_2-g(\na_XY,w_3)w_3.$$
 Thus
$$\theta_*\Tilde A_{1X}Y=
-g'(\na'_{\theta_*X}\theta_*Y,{\theta_*w_1}){\theta_*w_1}
-g'(\na'_{\theta_*X}\theta_*Y,{\theta_*w_2}){\theta_*w_2}=A_{1\theta_*X}\theta_*Y$$
for $\tilde\pi_1$-basic vector fields $X$ and $Y$,
where $g'$ denotes the metric on $\mathbb
CH^{2n+1}_{2s+1}$ and $\na'$ is the Levi-Civita connection of $g'$.

Let $X$ be the $\tilde\pi_1$-horizontal lifting along the fibre
$\tilde\pi_1^{-1}(\tilde\pi_1(p))$
of some unit vector in $T_{\tilde\pi_1(p)}\mathbb HH^n_s$.
Let $Y_1$, $Y_2$, $Y_3$ be the $\tilde\pi_1$-horizontal liftings along the fibre
$\tilde\pi_1^{-1}(\tilde\pi_1(p))$ of
$\tilde\pi_{1*}\tilde A_{1X}w_1$, $\tilde\pi_{1*}\tilde A_{1X}w_2$,
$\tilde\pi_{1*}\tilde A_{1X}w_3$,
respectively.
Let $v_i=\tilde A_{1X}Y_i$ for $i\in\{1,2,3\}$. As in Theorem \ref{p:3},
we choose $w_3=g(X,X)^{-1}\left(\na_{v_1}v_2\right) (p)$,
which implies that $v_3=\na_{v_1}v_2$ (see Lemma \ref{l:C}).

We remark that $v_3=\tilde A_{1X}Y_3$ is a
$\theta$-vertical vector field along the fibre
$\theta^{-1}(\theta(p))$. Indeed, we have
\begin{eqnarray*}
  \theta_*\left(\tilde A_{1X}Y_3(p')\right)&=&
      \left(A_{1\theta_*X}\theta_*Y_3\right)(\theta(p'))=
      \left(A_{1\theta_*X}\theta_*Y_3\right)(\theta(p))=
      \theta_*\left(\tilde A_{1X}Y_3(p)\right)\\
   &=&\theta_*(A_{1X}A_{1X}w_3)=g(X,X)\theta_*w_3=0
\end{eqnarray*}
for any $p'\in\theta^{-1}(\theta(p))$.

Since $v_1$, $v_2$ are orthogonal to the vertical vector field $v_3$ along
$\theta^{-1}(\theta(p))$, we see that $v_1$, $v_2$ are $\theta$-horizontal.
Since $\theta_*\left(\tilde A_{1X}Y_i(p')\right)=
\left(A_{1\theta_*X}\theta_*Y_i\right)(\theta(p'))$
for $p'\in\theta^{-1}(\theta(p))$ and for $i\in\{1,2\}$,
we obtain that $v_1$, $v_2$ are
$\theta$-basic vector fields along $\theta^{-1}(\theta(p))$.
Thus
$h\na_{v_3}v_1=A_{v_1}v_3$ along $\theta^{-1}(\theta(p))$.
Here $h$ and $v$ denote the $\theta$-horizontal
and $\theta$-vertical projections, respectively.
We also obtain that  $v\na_{v_3}v_1=-g(\na_{v_3}v_1,v_3)v_3=0$.
Therefore, $A_{v_1}v_3=\na_{v_3}v_1=v_2$ along $\theta^{-1}(\theta(p))$.

We shall prove that $\Tilde A_{1X}v_3=A_Xv_3$ along $\theta^{-1}(\theta(p))$ for every
 $\tilde\pi_1$-basic vector field $X$ along $\tilde\pi_1^{-1}(\tilde\pi_1(p))$.
We first obtain along $\theta^{-1}(\theta(p))$ that
   $$\Tilde A_{1X}v_3=\na_Xv_3+g(\na_Xv_3,v_1)v_1+g(\na_Xv_3,v_2)v_2+g(\na_Xv_3,v_3)v_3,$$
\begin{eqnarray*}
   g(\na_Xv_3,v_1)&=&g(A_Xv_3,v_1)=-g(v_3,A_Xv_1)=
      g(v_3,A_{v_1}X)=-g(A_{v_1}v_3,X)\\
     &=& -g(v_2,X)=0
\end{eqnarray*}
for a $\tilde\pi_1$-basic vector field  $X$ along $\tilde\pi_1^{-1}(\tilde\pi_1(p))$.
Analogously, we get $g(\na_Xv_3,v_2)=0$. Thus
$$\Tilde A_{1X}v_3=\na_Xv_3+g(\na_Xv_3,v_3)v_3=A_Xv_3$$
along $\theta^{-1}(\theta(p))$ for every
 $\tilde\pi_1$-basic vector field $X$ along $\tilde\pi_1^{-1}(\tilde\pi_1(p))$.

Let $\tilde{\mathcal L}= \{L_0=X, \Tilde A_{1L_0}v_1, \Tilde
A_{1L_0}v_2, \Tilde A_{1L_0}v_3, \dots, L_{n-1}, \Tilde
A_{1L_{n-1}}v_1, \Tilde A_{1L_{n-1}}v_2, \Tilde A_{1L_{n-1}}v_3\}$
be an orthonormal basis of $\tilde{\mathcal H}_{1}$ along the fibre
$\tilde\pi_1^{-1}(\tilde\pi_1(p))$ constructed as in Theorem
\ref{p:3}, for the semi-Riemannian submersion $\tilde\pi_1$.
From the proof of Theorem \ref{p:3}, we have
$$g(\tilde A_{1\tilde A_{1L_j}v_1}v_3,\tilde A_{1L_l}v_2)=0$$
 for $j\not=l$, and
$$g(\tilde A_{1\tilde A_{1L_j}v_1}v_3,L_t)=0$$
 for $0\leq j,t\leq n-1$.
We then obtain along $\tilde\pi_1^{-1}(\tilde\pi_1(p))$ that
\begin{eqnarray*}
g(\tilde A_{1\tilde A_{1L_j}v_1}v_3,\tilde A_{1L_j}v_2)&=&
-g(v_3,\tilde A_{1\tilde A_{1L_j}v_1}\tilde A_{1L_j}v_2)\\
&=&-g(v_3,\na_{v_1}v_2)g(L_j,L_j)\\
&=&-g(v_3,v_3)g(L_j,L_j)=-g(v_2,v_2)g(L_j,L_j)\\
&=&g(\tilde A_{1L_j}v_2,\tilde A_{1L_j}v_2),
\end{eqnarray*}
from which follows $\tilde A_{1L_j}v_2=\tilde A_{1\tilde A_{1L_j}v_1}v_3$.
Hence $\tilde A_{1L_j}v_2=A_{\tilde A_{1L_j}v_1}v_3$, because
$\tilde A_{1L_j}v_1$ is $\tilde\pi_1$-basic.
We also have $\tilde A_{1L_j}v_3=A_{L_j}v_3$.

Let $\mathcal L=\tilde{\mathcal L}\cup\{v_1,v_2\}$. Summarizing
all the above, we obtain that
$$\mathcal L=\{L_0,A_{L_0}v_3, \tilde A_{1L_0}v_1, A_{\tilde A_{1L_0}v_1}v_3, \dots,
L_{n-1}, A_{L_{n-1}}v_3, \tilde A_{L_{n-1}}v_1,
A_{\tilde A_{1L_{n-1}}v_1}v_3, v_1, A_{v_1}v_3\}$$
is an orthonormal basis  of the
$\theta$-horizontal space $\mathcal H$ along the fibre $\theta^{-1}(\theta(p))$
and $\mathcal L$ satisfies
all conditions imposed in the construction of the basis
$\mathcal L$ in the proof of Theorem \ref{p:3}.
We notice that $v_3=A_XY_3$ along  $\theta^{-1}(\theta(p))$,
and that along   $\theta^{-1}(\theta(p))$, $Y_3$ is equal to the
$\theta$-horizontal lifting of $\theta_*A_Xw_3$.

Let $q\in H^{4n+3}_{4s+3}$.
Let
   $$\tilde{\mathcal L}'=\{L'_0, \Tilde A_{2L'_0}v'_1,
   \Tilde A_{2L'_0}v'_2, \Tilde A_{2L'_0}v'_3,
   \dots, L'_{n-1}, \Tilde A_{2L'_{n-1}}v'_1, \Tilde A_{2L'_{n-1}}v'_2,
   \Tilde  A_{2L'_{n-1}}v'_3\}$$
be an orthonormal basis of $\tilde{\mathcal H}_{2}$
along $\tilde\pi_2^{-1}(\tilde\pi_2(q))$ constructed  in the same
way as $\tilde{\mathcal L}$, but for the semi-Riemannian
submersion $\tilde\pi_2$ (see the proof of Theorem \ref{p:3}), in
such  a way that $g_p(L_\al,L_\bet)=g_q(L'_\al,L'_\bet)$ for
$0\leq\al,\bet\leq n-1$, $g_p(v_i,v_j)=g_q(v'_i,v'_j)$ for
$1\leq i,j\leq 3$, and  $v'_3(q)=\left(\na_{v'_1}v'_2\right)(q)$.
Let $\phi:T_pH^{4n+3}_{4s+3}\to T_qH^{4n+3}_{4s+3}$ be the linear map given by
$\phi(v_i)=v'_i$, $\phi(\tilde A_{1L_\al}v_i)=\tilde A_{2L'_\al}v'_i$ for
$0\leq\al\leq n-1$ and for $1\leq i\leq 3$.

By Corollary 2.3.14 in \cite{wol}, $\phi$ extends to an isometry
$f:H^{4n+3}_{4s+3}\to H^{4n+3}_{4s+3}$ such that $f(p)=q$ and $f_{*p}=\phi$.
By the proof of Theorem \ref{p:3}, we have
$f_*\tilde A_{1E}F=\tilde A_{2f_*E}f_*F$ for every $E$, $F\in T_pH^{4n+3}_{4s+3}$.
By the proof of Theorem \ref{t:eq} and by Theorem \ref{t:4}, $f$ induces an isometry on
$\mathbb CH^{2n+1}_{2s+1}$, denoted by
$\tilde f:\mathbb CH^{2n+1}_{2s+1}\to\mathbb CH^{2n+1}_{2s+1}$,
such that $\theta\circ f=\tilde f\circ\theta$.
Since the $\pi_1$-vertical space at $\theta(p)$ is spanned by
$\{\theta_*v_1,\theta_*v_2\}$, since the $\pi_2$-vertical space at $\theta(q)$ is
spanned by
$\{\theta_*v'_1,\theta_*v'_2\}$, and since $\tilde f_*(\theta_*v_i)=\theta_*v'_i$,
for $i\in\{1,2\}$,
we see that $\tilde f_*$ maps the $\pi_1$-vertical space at  $\theta(p)$
into the $\pi_2$-vertical space at $\theta(q)$.
For  $\tilde\pi_1$-horizontal vectors $X$ and $Y$ we obtain
\begin{eqnarray*}
\tilde f_*A_{1\theta_*X}\theta_*Y&=&\tilde f_*\theta_*\tilde A_{1X}Y
=\theta_*f_*\tilde A_{1X}Y\\
&=&\theta_*\tilde A_{2f_*X}f_*Y=A_{2\theta_*f_*X}\theta_*f_*Y\\
&=&A_{2\tilde f_*(\theta_*X)}\tilde f_*(\theta_*Y).
\end{eqnarray*}
Therefore, by Theorem \ref{t:4}, we see that $\pi_1$ and $\pi_2$ are equivalent.
\end{proof}

\begin{remark}
We notice that our equivalence theorems can be applied, in particular,
to Riemannian submersions  from a sphere with totally geodesic fibres of dimension
less than or equal to $3$,
and for Riemannian submersions with complex totally geodesic fibres from a
complex projective space. Unlike those in  \cite{esc}, \cite{esco}, \cite{ran},
our proofs of the equivalence theorems are intrinsic,
we do not need to assume the existence
of any specific structure on the base space, such as
complex or quaternionic one. In Theorem \ref{t:eq2}, we need to assume only
that the fibres are 2-dimensional and that the induced metrics on fibres are negative
definite.
\end{remark}

Summarizing all results above, we now prove the main theorems.
\begin{proof}[Proof of Theorem \rm{\ref{t:1.1}}]
If $s+r'>1$, then $H^{n+r}_{s+r'}$ is simply connected and hence,
by Theorem \ref{p:3},
$B$ is an isotropic semi-Riemannian manifold and $r\in\{1,3\}$.
By Propositions \ref{p:8} and \ref{p:9}, we see that the base
space of the semi-Riemannian submersion is
isometric to a complex pseudo-hyperbolic space
if the dimension of fibres is one, or
to a quaternionic pseudo-hyperbolic space if the dimension of fibres is $3$.
In Theorem \ref{t:eq} we solved the equivalence problem.
The existence problem is solved by
the explicit construction given in the preliminaries (see Examples 1 and 2).

If $s+r'=1$, then either (i) $s=1$, $r'=0$, or (ii) $s=0$, $r'=1$.
Since the fibres are assumed to be negative definite, (i) cannot
occur.

(ii) If  $s=0$, $r'=1$, then $\pi$ is a
semi-Riemannian submersion from an anti-de Sitter space
onto a Riemannian manifold.
By \cite{mag}, $\pi$ is equivalent to the canonical submersion
$\pi:H^{2m+1}_1\to\mathbb CH^m$. This falls in the case (a).
\end{proof}

\begin{proof}[Proof of Theorem \rm{\ref{t:1.2}}]
If the dimension of the fibres is less than or
equal to $3$, then, by Theorem \ref{t:1.1},
$\pi$ is equivalent to the canonical semi-Riemannian submersions:
\begin{itemize}
\item[(a)] $H^{2m+1}_{2t+1}\to\mathbb CH^m_t$,  $0\leq t\leq m,$ or
\item[(b)] $H^{4m+3}_{4t+3}\to\mathbb HH^m_t$, $0\leq t\leq m.$
\end{itemize}

Now we assume that the dimension of the fibres is greater than or equal to 4.

(A)\ \ \ \ If we assume that the dimension of the fibres is greater
than or equal to $4$ and $B$ is an isotropic semi-Riemannian
manifold with non-constant curvature, then, by Proposition
\ref{p:8}, $B$ is isometric to $\mathbb CaH^2_t$, $t\in\{0,1,2\}$,
and the dimension of the fibres is $r=r'=7$. By Proposition
\ref{cayley}, there are no such semi-Riemannian submersions with
base space $\mathbb CaH^2_t$. Therefore, the assumptions (A) and
$r\geq 4$ imply that $B$ has constant curvature, and hence, by
Proposition \ref{p:9}, we obtain
$s=\mathrm{index}(B)\in\{0,\dim(B)\}$.

(B)\ \ \ If $\mathrm{index}(B)=0$ and $r\geq 4$, then, by \cite{bad},
the semi-Riemannian submersion $\pi$
is equivalent to the canonical
semi-Riemannian submersion ${H^{15}_7\to H^8(-4)}$.
If $\mathrm{index}(B)=\dim(B)$, then, by Lemma
\ref{l:new}, we get $r'=r$.
By changing the signs of the metrics on the base and on the total
space, $\pi$ becomes a Riemannian submersion with
connected totally geodesic fibres from a sphere onto
a Riemannian manifold. So, by \cite{esc} and \cite{ran},
 one obtains the conclusion.
\end{proof}

\begin{proof}[Proof of Theorem \rm{\ref{t:1.3}}]
Let
$\theta:H^{2n+1}_{2s+1}\to\mathbb CH^n_s$
be the canonical semi-Riemannian submersion.
By Theorem 2.5 in \cite{esco}, one obtains that
$\pi\circ\theta:H^{2n+1}_{2s+1}\to B$ is a semi-Riemannian
submersion with connected totally geodesic fibres.

(A)\ \ \ If the dimension of the fibres of $\pi$ is $r$ and $1\leq r\leq 2$, then
the dimension of the fibres of the semi-Riemannian submersion $\pi\circ\theta$
is less than or equal to $3$ and greater than or equal to $2$.
By Theorem \ref{t:1.1}, $B$ is isometric
to $\mathbb HH^m_t$ and $2n+1=4m+3$, $2s+1=4t+3$. Then $n=2m+1$, $s=2t+1$.
By Theorem \ref{t:eq2}, we see that
$\pi:\mathbb CH^{2m+1}_{2t+1}\to\mathbb HH^m_t$
is equivalent to the canonical semi-Riemannian
submersion.

(B) and (C)\ \ \ If $B$ is an isotropic semi-Riemannian
manifold or if $\mathrm{index}(B)\in\{0, \dim B\}$, then, by Theorem
\ref{t:1.2}, $\pi\circ\theta$ is equivalent
to one of the following canonical semi-Riemannian submersions:
\begin{itemize}
\item[] $H^{2m+1}_{2t+1}\to\mathbb CH^m_t$,\ \  $0\leq t\leq m$;
\item[] $H^{4m+3}_{4t+3}\to\mathbb HH^m_t$,\ \  $0\leq t\leq m$;
\item[] $H^{15}_{7+8t}\to H^8_{8t}(-4)$,\ \  $t\in\{0,1\}$.
\end{itemize}
If the dimension of the fibres of $\pi$ is greater than or equal to
$3$, then the dimension of the fibres of $\pi\circ\theta$ is greater
than or equal to $4$.
 Hence, in this case, $\pi\circ\theta$ is
equivalent to $H^{15}_{7+8t}\to H^8_{8t}(-4)$, $t\in\{0,1\}$.
 For $t=1$, the semi-Riemannian submersion $\pi$ is, after a change of
signs of the metrics on the total space and on the base space, of type $\pi:\mathbb
CP^7\to S^8(4)$. For $t=0$,  $\pi$ is of type $\pi:\mathbb
CH^7_3\to H^8(-4)$. In \cite{ran} (for case t=1) and \cite{bad}
(for case t=0), it is proved that there
are no such semi-Riemannian submersions with totally geodesic fibres.
 We proved that the dimension of fibres of $\pi$ is less than or equal to $2$.
\end{proof}

\begin{proof}[Proof of Theorem \rm{\ref{t:1.4}}]
We suppose that there are such semi-Riemannian submersions.
It is well-known that any quaternionic
submanifold in $\mathbb HH^n_s$  is totally geodesic.
Let $\eta:H^{4n+3}_{4s+3}\to\mathbb HH^n_s$,
$\xi:\mathbb CH^{2n+1}_{2s+1}\to\mathbb HH^n_s$,
be the canonical semi-Riemannian submersions.
By Theorem 2.5 in \cite{esco}, we see that
$\pi\circ\eta:H^{4n+3}_{4s+3}\to B$ is a semi-Riemannian
submersion with connected totally geodesic fibres.
We remark that the dimension of the fibres of $\pi\circ\eta$
is greater than or equal to $4$. Thus, by Theorem \ref{t:1.2},
we see that  $\pi\circ\eta$ is
equivalent to the canonical semi-Riemannian submersion
$$ H^{15}_7\to H^8(-4),\ \mathrm{or}\  H^{15}_{15}\to H^8_8(-4).$$
It follows that $\pi$ is one of the following types:
\begin{itemize}
  \item[(i)] $\pi:\mathbb HH^3_1\to H^8(-4)$, or
  \item[(ii)] $\pi:\mathbb HH^3_3\to H^8_8(-4)$.
\end{itemize}
In \cite{ucc}, Ucci proved that there are no
Riemannian submersions with fibres $\mathbb HP^1$
from $\mathbb HP^3$ onto $S^8(4)$. Therefore, Case (ii) is not possible.

The fibres of semi-Riemannian submersion $\pi\circ\xi:\mathbb CH^{7}_{3}\to H^8(-4)$
are totally geodesic by Theorem 2.5 in \cite{esco}, and complex submanifolds,
since the horizontal
lifting of the tangent space of the quaternionic line
$\pi^{-1}(\pi(p))$ is invariant under the canonical complex
structure on $\mathbb CH^7_3$.  By \cite{bad},
there are no semi-Riemannian submersions with complex totally
geodesic fibres from  $\mathbb CH^7_3$ onto $H^8(-4)$. Thus Case (i) is impossible.
\end{proof}

\end{document}